\documentclass[12pt]{article}
\usepackage{amssymb, amsmath, amsthm, amscd}
\usepackage[pdftex]{graphics}
\usepackage[utf8]{inputenc}
\usepackage[all,cmtip]{xy}
\usepackage{enumitem}
\usepackage{setspace}
\usepackage{mathdots}

\usepackage[mathscr]{euscript}
\usepackage[colorlinks=true]{hyperref}
\hypersetup{colorlinks   = true}
\hypersetup{linkcolor=blue}
\usepackage{tikz}

\begin{document}

\newtheorem{The}{Theorem}[section]
\newtheorem{Lem}[The]{Lemma}
\newtheorem{Prop}[The]{Proposition}
\newtheorem{Cor}[The]{Corollary}
\newtheorem{Rem}[The]{Remark}
\newtheorem{Obs}[The]{Observation}
\newtheorem{SConj}[The]{Standard Conjecture}
\newtheorem{Titre}[The]{\!\!\!\! }
\newtheorem{Conj}[The]{Conjecture}
\newtheorem{Question}[The]{Question}
\newtheorem{Prob}[The]{Problem}
\newtheorem{Def}[The]{Definition}
\newtheorem{Not}[The]{Notation}
\newtheorem{Claim}[The]{Claim}
\newtheorem{Conc}[The]{Conclusion}
\newtheorem{Ex}[The]{Example}
\newtheorem{Fact}[The]{Fact}
\newtheorem{Formula}[The]{Formula}
\newtheorem{Formulae}[The]{Formulae}
\newtheorem{The-Def}[The]{Theorem and Definition}
\newtheorem{Prop-Def}[The]{Proposition and Definition}
\newtheorem{Lem-Def}[The]{Lemma and Definition}
\newtheorem{Cor-Def}[The]{Corollary and Definition}
\newtheorem{Conc-Def}[The]{Conclusion and Definition}
\newtheorem{Terminology}[The]{Note on terminology}
\newcommand{\C}{\mathbb{C}}
\newcommand{\R}{\mathbb{R}}
\newcommand{\N}{\mathbb{N}}
\newcommand{\Z}{\mathbb{Z}}
\newcommand{\Q}{\mathbb{Q}}
\newcommand{\Proj}{\mathbb{P}}
\newcommand{\Rc}{\mathcal{R}}
\newcommand{\Oc}{\mathcal{O}}
\newcommand{\Vc}{\mathcal{V}}
\newcommand{\Id}{\operatorname{Id}}
\newcommand{\pr}{\operatorname{pr}}
\newcommand{\rk}{\operatorname{rk}}
\newcommand{\del}{\partial}
\newcommand{\delbar}{\bar{\partial}}
\newcommand{\Cdot}{{\raisebox{-0.7ex}[0pt][0pt]{\scalebox{2.0}{$\cdot$}}}}
\newcommand\nilm{\Gamma\backslash G}
\newcommand\frg{{\mathfrak g}}
\newcommand{\fg}{\mathfrak g}
\newcommand{\Oh}{\mathcal{O}}
\newcommand{\Kur}{\operatorname{Kur}}
\newcommand\gc{\frg_\mathbb{C}}

\begin{center}

{\Large\bf A Twisted Adiabatic Limit Approach to Vanishing Theorems for Complex Line Bundles}

\end{center}

\begin{center}

{\large Dan Popovici}

\end{center}

\vspace{1ex}

\noindent{\small{\bf Abstract.} Given an $n$-dimensional compact complex Hermitian manifold $X$, a $C^\infty$ complex line bundle $L$ equipped with a connection $D$ whose $(0,\,1)$-component $D''$ squares to zero and a real-valued function $\eta$ on $X$, we prove that the $D''$-cohomology group of $L$ of any bidegree $(p,\,q)$ such that either $(p>q \hspace{1ex}\mbox{and}\hspace{1ex} p+q\geq n+1)$ or $(p<q \hspace{1ex}\mbox{and}\hspace{1ex} p+q\leq n-1)$ vanishes when two extra hypotheses are made. The first hypothesis requires a certain real-valued, not necessarily closed, $(1,\,1)$-form depending on $p,\,q$, on the curvature of $D$ and on a $(1,\,1)$-form induced by $\eta$ to be positive definite. The second hypothesis requires the norm of $\partial\eta$ to be small relative to $|\eta|$. This theorem, for which we also give a number of variants, is proved by generalising our very recent twisted adiabatic limit construction for complex structures to connections on complex line bundles. This twisting of $D$ induces first-order differential operators acting on the $L$-valued forms, for which we obtain commutation relations involving their formal adjoints, and two twisted Laplacians for which we obtain a comparison formula reminiscent of the classical Bochner-Kodaira-Nakano identity. The main features of our results are that $X$ need not be K\"ahler, $L$ need not be holomorphic and the types of $C^\infty$ functions that $X$ supports play a key role in our hypotheses, thus capturing some of their links with the geometry of manifolds.}

\vspace{1ex}

\section{Introduction}\label{secion:introduction}

It is a standard fact that for a holomorphic line bundle $L$ on a compact complex manifold $X$, the ampleness of $L$ is equivalent to the existence of a $C^\infty$ Hermitian fibre metric $h$ whose curvature form $i\Theta_h(L)$ ($=$ the curvature form of the Chern connection $D$ of $(L,\,h)$) is positive definite everywhere on $X$. Moreover, the celebrated Kodaira Embedding Theorem stipulates that $X$ is projective if and only if an ample line bundle $L$ exists thereon.

There is a huge literature about vanishing theorems on projective and more general manifolds. For the reader's convenience, we recall a few milestones.

\vspace{2ex}

\noindent {\bf Classical vanishing theorems} {\it Let $L$ be a holomorphic line bundle on a complex {\bf projective} manifold $X$ with $\mbox{dim}_\C X =n$.

\vspace{1ex}

(i)\, ({\bf Kodaira's vanishing theorem [Kod53]}) If $L$ is {\bf ample}, then $H^{n,\,q}(X,\,L) = \{0\}$ whenever $q\geq 1$.
  
\vspace{1ex}

(ii)\, ({\bf Akizuki-Nakano vanishing theorem [AN54]}) If $L$ is {\bf ample}, then $H^{p,\,q}(X,\,L) = \{0\}$ whenever $p+q\geq n+1$.

\vspace{1ex}

(iii)\, ({\bf Grauert-Riemenschneider vanishing theorem [GR70]}) If $L$ is {\bf semi-ample} and {\bf big}, then $H^{n,\,q}(X,\,L) = \{0\}$ whenever $q\geq 1$.}

\vspace{2ex}

In particular, the Akizuki-Nakano and the Grauert-Riemenschneider vanishing theorems generalise Kodaira's theorem in different ways: extra bidegrees $(p,\,q)$ besides $(n,\,q)$ are captured in the former case; a weaker positivity assumption than ampleness is made on $L$ in the latter case.

Actually, the original form of the Grauert-Riemenschneider vanishing theorem [GR70, Satz 2.1.] ensures the vanishing of $H^{n,\,q}(X,\,L)$ whenever $X$ is merely {\bf Moishezon} (and it need not even be a manifold, it suffices for it to be an irreducible compact complex space) and $L$ is {\bf quasi-positive} (in the sense that it admits a $C^\infty$ Hermitian fibre metric $h$ whose curvature form $i\Theta_h(L)$ is non-negative everywhere on $X$ and positive definite at some point). More generally, $L$ need not even be a line bundle, it can merely be a quasi-positive torsion-free coherent analytic sheaf of rank $1$ on $X$.

Grauert and Riemenschneider went on to ask ([GR70, p. 277]), in a question that came to be known as the {\it Grauert-Riemenschneider Conjecture}, whether the Moishezon assumption on $X$ is superfluous, namely whether it is implied by the existence of a quasi-positive torsion-free coherent analytic sheaf of rank $1$ on $X$. (Such a sheaf can always be assumed locally free of rank $1$ after desingularisation.) This conjecture was answered in the affirmative by Siu in [Siu84] and [Siu85]. Demailly went on to reprove it as a corollary of a stronger result ([Dem85, th\'eor\`eme 0.8.]), obtained in turn as a consequence of his {\it Morse inequalities}, that further weakens the quasi-positivity assumption on the Hermitian holomorphic line bundle $(L,\,h)$ which is required to exist on $X$ to the condition \begin{eqnarray}\label{eqn:Demailly_G-R_hypothesis}\int\limits_{X(L,\,\leq 1)}\bigg(i\Theta_h(L)\bigg)^n > 0,\end{eqnarray} where $X(L,\,\leq 1)$ is the set of points in $X$ at which the curvature form $i\Theta_h(L)$ has at most one negative eigenvalue.  

\vspace{1ex}

The literature also abounds in generalisations of these results to the case of holomorphic vector bundles $F$ of arbitrary rank $r\geq 1$, such as:

\vspace{1ex}

-the {\it Nakano vanishing theorem} of [Nak55], giving the same conclusion as Kodaira's vanishing theorem in the case where $F$ (replacing $L$) is supposed {\it Nakano positive};

\vspace{1ex}

-the {\it Le Potier vanishing theorem} of [LP75], giving a suitably modified conclusion (i.e. $H^{p,\,q}(X,\,F) = \{0\}$ whenever $p+q\geq n+r$) from the Akizuki-Nakano vanishing theorem in the case where $F$ is supposed {\it Griffiths positive}.     

\vspace{2ex}

A key feature of the outstanding results of Siu in [Siu84] and [Siu85] and Demailly in [Dem85] is that the ambient manifold $X$ is not assumed K\"ahler, let alone projective. Its Moishezon property (which does not imply K\"ahlerianity) is only deduced a posteriori if an appropriate positivity assumption is made on $L$. However, no such assumption is made in Demailly's {\it holomorphic Morse inequalities} of [Dem85, th\'eor\`eme 0.1.]; they hold on arbitrary compact complex manifolds $X$.

\vspace{1ex}

Consistent with the use of his holomorphic Morse inequalities for a holomorphic line bundle $L$ (hence for an integral cohomology class $c_1(L)$ of type $(1,\,1)$, the first Chern class of $L$) in getting a sufficient criterion for $X$ to be Moishezon, Demailly proposed a far-reaching conjecture in the general context of possibly transcendental cohomology classes.

\begin{Conj}(Demailly's transcendental Morse inequalities conjecture)\label{Conj:Demailly_transcendental-Morse} Let $X$ be an $n$-dimensional compact complex manifold. Suppose there exists a real $d$-closed $C^\infty$ $(1,\,1)$-form $\alpha$
  on $X$ such that \begin{eqnarray}\label{eqn:Demailly_tr-Morse_hypothesis}\int\limits_{X(\alpha,\,\leq 1)}\alpha^n > 0,\end{eqnarray} where $X(\alpha,\,\leq 1)$ is the set of points in $X$ at which $\alpha$ has at most one negative eigenvalue.

  Then, the Bott-Chern cohomology class $\{\alpha\}_{BC}$ contains a {\bf K\"ahler current}. In particular, $X$ is a {\bf class ${\cal C}$ manifold}.

\end{Conj}

\vspace{2ex}

Recall that a real current $T$ of bidegree $(1,\,1)$ on $X$ is said to be a {\it K\"ahler current} if $dT=0$ and if, for a given Hermitian metric $\omega$ on $X$ (seen as a $C^\infty$ positive definite $(1,\,1)$-form on $X$ -- such an $\omega$ always exists), there exists a constant $\varepsilon>0$ such that $T\geq\varepsilon\omega$. Further recall that, by Theorem 3.4. in [DP04], the existence of such a current on $X$ is equivalent to $X$ being a {\it class ${\cal C}$ manifold}, namely to $X$ being bimeromorphically equivalent to a compact K\"ahler manifold.

\vspace{2ex}

This paper is motivated by an attempt to make headway towards an eventual resolution of Conjecture \ref{Conj:Demailly_transcendental-Morse}. One of its main results, reminiscent of the Akizuki-Nakano vanishing theorem [AN54] and of the Grauert-Riemenschneider vanishing theorem [GR70], but with no assumption on $X$ and with a different kind of positivity assumption on a possibly non-holomorphic line bundle $L$, involving smooth functions on $X$ (a key aspect of our approach), is the following consequence of Theorem \ref{The:vanishing_cohomology_bundle_compact}.

\begin{The}\label{The:vanishing_cohomology_bundle_compact_introd} Let $X$ be a {\bf compact} complex manifold with $\mbox{dim}_\C X=n\geq 2$ and let $L$ be a $C^\infty$ complex line bundle on $X$ equipped with a $C^\infty$ Hermitian fibre metric $h$ and with an {\bf h-compatible} linear connection $D=D' + D''$ such that $D''^2=0$. Fix integers $p,q\in\{0,\dots , n\}$ such that either $(p>q \hspace{1ex}\mbox{and}\hspace{1ex} p+q\geq n+1)$ or $(p<q \hspace{1ex}\mbox{and}\hspace{1ex} p+q\leq n-1)$.

  Suppose there exists a family $(\eta_\varepsilon)_\varepsilon$ of non-vanishing $C^\infty$ functions $\eta_\varepsilon:X\longrightarrow\R$ such that: 

\vspace{1ex}

(i)\, the $C^\infty$ $(1,\,1)$-form $$\omega_\varepsilon = \omega_{p,\,q,\,\eta_\varepsilon}:= \frac{1}{p-q}\,i\Theta(D)^{1,\,1} + \gamma_{\eta_\varepsilon}$$ is {\bf positive definite} at every point of $X$ for every $\varepsilon>0$, where $\Theta(D)^{1,\,1}$ is the part of type $(1,\,1)$ of the curvature form $\Theta(D)$ of $D$ and $\gamma_{\eta_\varepsilon}:=\frac{2}{\eta_\varepsilon^2}\,i\partial\eta_\varepsilon\wedge\bar\partial\eta_\varepsilon - \frac{1}{\eta_\varepsilon}\,i\partial\bar\partial\eta_\varepsilon$;

\vspace{1ex}

(ii)\, the pointwise $\omega_\varepsilon$-norm $|\partial\eta_\varepsilon| = |\partial\eta_\varepsilon|_{\omega_\varepsilon}$ of the $(1,\,0)$-form $\partial\eta_\varepsilon$ is {\bf asymptotically small} relative to $|\eta_\varepsilon|$ in the sense that $C_1(\eta_\varepsilon):=\sup\limits_X\frac{|\partial\eta_\varepsilon|}{|\eta_\varepsilon|}$ converges to $0$ as $\varepsilon\downarrow 0$;

\vspace{1ex}

(iii)\, there exists a constant $C>0$ such that $C(\varphi_\varepsilon):=\sup\limits_X|i\partial\bar\partial\varphi_\varepsilon|_{\omega_\varepsilon}\leq C$ for every $\varepsilon>0$, where $\varphi_\varepsilon:=-\log|\eta_\varepsilon|$.

 \vspace{1ex}

Then, the $D''$-cohomology group of $L$ of bidegree $(p,\,q)$ {\bf vanishes}, namely $H^{p,\,q}_{D''}(X,\,L) = \{0\}$.

\end{The}
  
By the connection $D$ being {\bf h-compatible} we mean that it is a {\it Hermitian} connection, namely compatible with the Hermitian metric $h$ in the sense that it satisfies the standard condition (\ref{eqn:Hermitian-connection_def}). For an h-compatible connection, the integrability of its $(0,\,1)$-part $D''$ (i.e. the property $D''^2=0$) implies the integrability of its $(1,\,0)$-part $D'$ (i.e. $D'^2 = 0$), hence it implies that the curvature form $\Theta(D)$ is reduced to its component $\Theta(D)^{1,\,1}$ of type $(1,\,1)$. In particular, $\Theta(D)^{1,\,1}$ is closed in this case. Meanwhile, the real $(1,\,1)$-form $\gamma_{\eta_\varepsilon}$ need not be closed, although $(1/\eta_\varepsilon)\,\gamma_{\eta_\varepsilon}$ is closed. In fact, an immediate computation yields: \begin{eqnarray}\label{eqn:ddbar_1-over-eta}i\partial\bar\partial\bigg(\frac{1}{\eta}\bigg) = \frac{1}{\eta}\,\bigg(\frac{2}{\eta^2}\,i\partial\eta\wedge\bar\partial\eta - \frac{1}{\eta}\,i\partial\bar\partial\eta\bigg) =  \frac{1}{\eta}\,\gamma_\eta\end{eqnarray} for any non-vanishing $C^\infty$ function $\eta_\varepsilon:X\longrightarrow\R$.

  This discussion shows that neither the forms $\omega_\varepsilon$ of Theorem \ref{The:vanishing_cohomology_bundle_compact_introd}, nor the forms $(1/\eta_\varepsilon)\,\omega_\varepsilon$, need be closed, so they need not define K\"ahler metrics on $X$. Thus, one of the key features of Theorem \ref{The:vanishing_cohomology_bundle_compact_introd} is that it constitutes a vanishing result on possibly non-K\"ahler compact complex manifolds in the spirit of [Siu84], [Siu85] and [Dem85]. However, unlike in some of these references, even a posteriori the manifolds involved in our results need not be either K\"ahler or Moishezon.

  This is made possible by the positivity assumption bearing not directly on the curvature form $i\Theta(D)^{1,\,1}$, but on the (possibly non-closed) combinations $\omega_\varepsilon$ of it with the forms $\gamma_{\eta_\varepsilon}$ induced by the functions $\eta_\varepsilon$. Thus, the vanishing of certain cohomology groups of certain complex line bundles $L\longrightarrow X$ arises as a partial consequence of the types of smooth functions that the compact complex manifold $X$ supports. The study of the relationships between these functions and the complex line bundles on (hence the geometry of) complex manifolds is one of the central planks of our approach. In our opinion, it deserves further probing.

  Another peculiarity of our results is that they deal with complex line bundles $L$ that are not assumed holomorphic, but only $C^\infty$. One of the motivations for this generality is the need, in some cases, to deal with {\it transcendental} cohomology classes $\{\alpha\}$ of type $(1,\,1)$ as those involved in Conjecture \ref{Conj:Demailly_transcendental-Morse} or in Siu's conjecture on the invariance of the plurigenera in K\"ahler families of compact complex manifolds. In that case, positive integer multiples $k\,\{\alpha\}$ can be approximated by {\it integral} cohomology classes $\{\alpha_k\}$ that may no longer be of type $(1,\,1)$. One then obtains a sequence of $C^\infty$ but, in general, non-holomorphic (although approximately holomorphic) complex line bundles $L_k$ whose first Chern classes are the $\{\alpha_k\}$'s, as explained in e.g. [Pop23].   

\vspace{2ex}

The proofs of Theorem \ref{The:vanishing_cohomology_bundle_compact_introd} and of the more general Theorem \ref{The:vanishing_cohomology_bundle_compact} proceed through the introduction of what we call the {\it twisted adiabatic limit} construction. The version thereof proposed in this paper for connections on complex vector bundles generalises the one introduced very recently in [Pop24] for complex structures and scalar-valued differential forms. Both constructions, new to our knowledge, use a possibly non-constant function $\eta$ to twist a connection $D$ or the Poincar\'e operator $d$, as opposed to a constant $h$ used earlier in [Pop19] and even earlier by various authors (see background discussions in [Pop19] and [Pop24]) in the different context of Riemannian foliations.

We now briefly outline the main results of this paper. The main construction is performed in $\S$\ref{section:main-def}. In the fairly general setting of a $C^\infty$ complex line bundle $L$ equipped with a connection $D = D' + D''$ (split into its components of types $(1,\,0)$ and $(0,\,1)$) over an $n$-dimensional complex manifold $X$, we associate with any non-vanishing $\C$-valued $C^\infty$ function $\eta$ on $X$ the twisted operators $D_\eta^{1,\,0} = \theta_\eta D'\theta_\eta^{-1}:C^\infty_{p,\,q}(X,\,L)\longrightarrow C^\infty_{p+1,\,q}(X,\,L)$ and $D_\eta^{0,\,1} = \theta_\eta D''\theta_\eta^{-1}:C^\infty_{p,\,q}(X,\,L)\longrightarrow C^\infty_{p,\,q+1}(X,\,L)$ acting on $L$-valued $C^\infty$ forms of any bidegree $(p,\,q)$, where $\theta_\eta$ is the operator acting pointwise on these forms by multiplication by $\eta^p$.

Once $L$ has been given a Hermitian fibre metric $h$ and $X$ a Hermitian metric $\omega$, the induced $L^2$-inner product determines formal adjoints of these operators. These induce Laplace-type elliptic differential operators of order two $\Delta''_\eta:=[D_\eta^{0,\,1},\,(D_\eta^{0,\,1})^\star]$ and $\Delta'_\eta:=[\overline{D_\eta^{0,\,1}},\,\overline{D_\eta^{0,\,1}}^{\,\star}]$ acting on every space $C^\infty_{p,\,q}(X,\,L)$ with values in itself. When the $(0,\,1)$-component $D''$ of the original connection $D$ is {\it integrable} (i.e. $D''^2 = 0$), so is the twisted operator $D^{0,\,1}_\eta$. Hence, it defines finite-dimensional cohomology groups $H^{p,\,q}_{D^{0,\,1}_\eta}(X,\,L)$ that turn out to be isomorphic to the $D''$-cohomology groups $H^{p,\,q}_{D''}(X,\,L)$ and, via Hodge-type isomorphisms, to the kernels ${\cal H}_{\Delta''_\eta}^{p,\,q}(X,\,L)$ of $\Delta''_\eta$ of the same bidegree. (See Observation \ref{Obs:D''_cohom_Hodge}.) The takeaway from this twisted adiabatic limit construction (whose ``limiting'' feature manifests itself when the twisting function $\eta$ is replaced by a family $(\eta_\varepsilon)_{\varepsilon}$ of functions and $\varepsilon$ is let to converge to some limit, possibly $0$, as in Theorem \ref{The:vanishing_cohomology_bundle_compact_introd}) is that we can study the $D''$-cohomology of $L$ by studying the $\eta$-twisted Laplacian $\Delta''_\eta$.

\vspace{1ex}

The study of $D_\eta^{0,\,1}$ and $\Delta''_\eta$ is carried on in $\S$\ref{section:twisted-commutation_D-eta_0-1} where {\it $\eta$-twisted commutation relations} are obtained. They express the adjoints of $D_\eta^{0,\,1}$ and its conjugate in terms of the anti-commutators between, on the one hand, $\Lambda = \Lambda_\omega$ (the adjoint of the multiplication of differential forms by the metric $\omega$), and, on the other hand, the conjugate of $D_\eta^{0,\,1}$, resp. $D_\eta^{0,\,1}$ itself. This is done in Proposition \ref{Prop:eta-twisted_commutation-relations_1-0_0-1}.

The {\it $\eta$-twisted commutation relations} lead to an {\it $\eta$-twisted Bochner-Kodaira-Nakano-type} identity (see Proposition \ref{Prop:eta-BKN_rough}) comparing the twisted Laplacians $\Delta''_\eta$ and $\Delta'_\eta$. This identity contains the curvature-type operator $i\,\bigg[[D^{0,\,1}_\eta,\,\overline{D^{0,\,1}_\eta}],\,\Lambda\bigg]$ which, perhaps surprisingly, is of order $1$ (unlike the order $0$ of the classical untwisted case), as shown by the computation contained in Proposition \ref{Prop:curvature_computation} and its proof. However, the curvature-type operator is seen to contain the zero-th order operator $[\Omega_{p,\,q}\wedge\cdot\,,\,\Lambda]$ when acting in bidegree $(p,\,q)$, where $\Omega_{p,\,q}$ is a $C^\infty$ $\C$-valued $(1,\,1)$-form on $X$ which, in the case where the function $\eta$ is {\it real-valued}, is given (cf. Observation \ref{Obs:eta_real-valued_curvature-form}) by the formula $$\Omega_{p,\,q}:= i\Theta(D)^{1,\,1} + (p-q)\,\gamma_\eta.$$ This formula features both the $(1,\,1)$-part of the curvature form of the connection $D$ and the real $(1,\,1)$-form $\gamma_\eta$ of (\ref{eqn:ddbar_1-over-eta}).

In $\S$\ref{section:vanishing_cohomology_bundle_compact}, as applications of our twisted adiabatic limit construction and our {\it $\eta$-twisted Bochner-Kodaira-Nakano-type} identity of the previous sections, we get vanishing theorems. Chief among them is Theorem \ref{The:vanishing_cohomology_bundle_compact}. We then deduce a number of variants of it: Corollaries \ref{Cor:vanishing_cohomology_bundle_compact} and \ref{Cor:vanishing_cohomology_bundle_compact_simpler}, as well as Theorem \ref{The:vanishing_cohomology_bundle_compact_introd}. In the spirit of Problem 1.1. proposed in [Pop24] and of Proposition 1.3. of [Pop24], we further give in Corollary \ref{Cor:spectrum-lower-bound}, as an application of our $\eta$-twisted Bochner-Kodaira-Nakano-type identity, a lower bound for the {\it positive} part of the spectrum of the $\eta$-twisted Laplacian $\Delta''_\eta$ in the bidegree $(n,\,0)$. This bidegree does not come within the scope of our other results such as Theorem \ref{The:vanishing_cohomology_bundle_compact}. In all these applications, the twisting function $\eta$ is supposed real-valued since it defines quantities involved in positivity hypotheses.

\section{Main definitions}\label{section:main-def}

Let $X$ be a complex manifold with $\mbox{dim}_\C X = n$ on which a Hermitian metric $\omega$ has been fixed. Let $L$ be a $C^\infty$ complex line bundle on $X$ equipped with a $C^\infty$ Hermitian fibre metric $h$ and with a linear connection $D = D' + D''$, where $D'$ and $D''$ are the components of types $(1,\,0)$, resp. $(0,\,1)$, of $D$. We do not assume $L$ to be holomorphic or $D$ to be compatible with $h$ unless otherwise specified.

Let $\eta:X\longrightarrow\C$ be a non-vanishing $C^\infty$ function on $X$. In this paper, we generalise to the operators $$D':C^\infty_{p,\,q}(X,\,L)\longrightarrow C^\infty_{p+1,\,q}(X,\,L) \hspace{3ex}\mbox{and}\hspace{3ex} D'':C^\infty_{p,\,q}(X,\,L)\longrightarrow C^\infty_{p,\,q+1}(X,\,L)$$ the construction performed in [Pop24] for their counterparts $\partial:C^\infty_{p,\,q}(X,\,\C)\longrightarrow C^\infty_{p+1,\,q}(X,\,\C)$ and $\bar\partial:C^\infty_{p,\,q}(X,\,\C)\longrightarrow C^\infty_{p,\,q+1}(X,\,\C)$ acting on scalar-valued forms on $X$.

\begin{Def}\label{Def:theta-D_eta} Fix an arbitrary degree $k\in\{0,\dots , 2n\}$.

  \vspace{1ex}

 (i)\, Let $\theta_\eta:\Lambda^k T^\star X\otimes L \longrightarrow \Lambda^k T^\star X\otimes L$ be the bijective linear operator defined pointwise on $L$-valued $k$-forms on $X$ by \begin{eqnarray}\label{eqn:theta_eta_def}\Lambda^k T^\star X\otimes L \ni u=\sum\limits_{p+q=k}u^{p,\,q} \longmapsto  \sum\limits_{p+q=k}\theta_\eta(u^{p,\,q}):=\sum\limits_{p+q=k}\eta^p\,u^{p,\,q}\in \Lambda^k T^\star X\otimes L.\end{eqnarray}

  (ii)\, For any degree $k\in\{0,\dots , 2n\}$, let $D_\eta:C^\infty_k(X,\,L)\longrightarrow C^\infty_{k+1}(X,\,L)$ be the first-order differential operator defined on the $L$-valued $C^\infty$ $k$-forms on $X$ by \begin{eqnarray}\label{eqn:D_eta}D_\eta = \theta_\eta D\theta_\eta^{-1}.\end{eqnarray}

\end{Def}  

This is the analogue in this vector bundle context of (3) and Definition 2.2. of [Pop24]. After splitting $D_\eta = D_\eta^{1,\,0} + D_\eta^{0,\,1}$ into its components of types $(1,\,0)$ and $(0,\,1)$, we get first-order linear differential operators: \begin{eqnarray}\label{eqn:D_eta_1-0_def}D_\eta^{1,\,0} = \theta_\eta D'\theta_\eta^{-1} = \eta D' - p\,\partial\eta\wedge\cdot\,:C^\infty_{p,\,q}(X,\,L)\longrightarrow C^\infty_{p+1,\,q}(X,\,L)\end{eqnarray} and \begin{eqnarray}\label{eqn:D_eta_0-1_def}D_\eta^{0,\,1} = \theta_\eta D''\theta_\eta^{-1} = D'' - \frac{p}{\eta}\,\bar\partial\eta\wedge\cdot\,:C^\infty_{p,\,q}(X,\,L)\longrightarrow C^\infty_{p,\,q+1}(X,\,L)\end{eqnarray} in every bidegree $(p,\,q)$. These explicit formulae are obtained in the same way as their scalar analogues in [Pop24, Proposition and Definition 2.4.].

The formal adjoints of these operators w.r.t. the $L^2$-inner product induced on $\bigoplus_{p,\,q}C^\infty_{p,\,q}(X,\,L)$ by the metrics $\omega$ of $X$ and $h$ of $L$ are \begin{eqnarray}\label{eqn:D_eta_1-0_star}(D_\eta^{1,\,0})^\star = D'^\star(\overline\eta\cdot\,) - p\,(\partial\eta\wedge\cdot\,)^\star:C^\infty_{p+1,\,q}(X,\,L)\longrightarrow C^\infty_{p,\,q}(X,\,L)\end{eqnarray} and \begin{eqnarray}\label{eqn:D_eta_0-1_star}(D_\eta^{0,\,1})^\star = D''^\star - \frac{p}{\overline\eta}\,(\bar\partial\eta\wedge\cdot\,)^\star:C^\infty_{p,\,q+1}(X,\,L)\longrightarrow C^\infty_{p,\,q}(X,\,L).\end{eqnarray}

As in the scalar case of [Pop24], we also consider the Laplacians: \begin{eqnarray}\label{eqn:Laplacians_D_eta_0-1_star}\nonumber\Delta''_\eta:=[D_\eta^{0,\,1},\,(D_\eta^{0,\,1})^\star]:C^\infty_{p,\,q}(X,\,L) & \longrightarrow & C^\infty_{p,\,q}(X,\,L) \\  \Delta'_\eta:=[\overline{D_\eta^{0,\,1}},\,\overline{D_\eta^{0,\,1}}^{\,\star}]:C^\infty_{p,\,q}(X,\,L) & \longrightarrow & C^\infty_{p,\,q}(X,\,L),\end{eqnarray} where the conjugate operator $\overline{D_\eta^{0,\,1}}:C^\infty_{p,\,q}(X,\,L)\longrightarrow C^\infty_{p+1,\,q}(X,\,L)$ is defined by requiring $\overline{D_\eta^{0,\,1}}u=\overline{D_{\eta,\,\overline{L}}^{0,\,1}\,\overline{u}}$ for every form $u\in C^\infty_{p,\,q}(X,\,L)$, while $\overline{u}\in C^\infty_{q,\,p}(X,\,\overline{L})$ is the conjugate of $u$ and $D_{\eta,\,\overline{L}}^{0,\,1}:C^\infty_{q,\,p}(X,\,\overline{L})\longrightarrow C^\infty_{q,\,p+1}(X,\,\overline{L})$ is defined in $\overline{L}$ analogously to $D_\eta^{0,\,1}$ in $L$. When no confusion is likely, we will simply denote $D_{\eta,\,\overline{L}}^{0,\,1}$ by $D_\eta^{0,\,1}$. It will be clear in context whether it is acting on an $L$-valued form or on an $\overline{L}$-valued one. 

A word of explanation is in order here. By $\overline{L}$ we mean the $C^\infty$ complex line bundle on $X$ that is the conjugate of $L$ in the sense that its transition functions are $\overline{g}_{\alpha\beta}$, where the $g_{\alpha\beta}$ are the transition functions of $L$ w.r.t. to a given collection $(U_\alpha,\,\tau_\alpha)_\alpha$ of local trivialisations. Now, the connection $D$ on $L$ induces the conjugate connection $D_{\overline{L}}= D_{\overline{L}}' + D_{\overline{L}}''$ on $\overline{L}$ whose components of type $(1,\,0)$, resp. $(0,\,1)$, are defined by $$D_{\overline{L}}'\overline{u}:=\overline{D''u} \hspace{3ex} \mbox{and} \hspace{3ex} D_{\overline{L}}''\overline{u}:=\overline{D'u}$$ for any $u\in C^{p,\,q}(X,\,L)$ and any bidegree $(p,\,q)$. (Each $u\in C^{p,\,q}(X,\,L)$ can be seen as a collection $(u_\alpha)_\alpha$ of $\C$-valued $(p,\,q)$-forms defined respectively on the open subsets $U_\alpha\subset X$ such that $u_\alpha = g_{\alpha\beta}\,u_\beta$ on $U_\alpha\cap U_\beta$ for all $\alpha,\,\beta$; $\overline{u}$ is then the $\overline{L}$-valued $(q,\,p)$-form defined by the collection $(\overline{u}_\alpha)_\alpha$ of locally defined $\C$-valued $(q,\,p)$-forms glued together by the relations $\overline{u}_\alpha = \overline{g}_{\alpha\beta}\,\overline{u}_\beta$ on $U_\alpha\cap U_\beta$ for all $\alpha,\,\beta$.)

It is by means of $D_{\overline{L}}'$ and $D_{\overline{L}}''$ that the conjugates of the operators $D^{1,\,0}_\eta:C^\infty_{p,\,q}(X,\,L)\longrightarrow C^\infty_{p+1,\,q}(X,\,L)$ and $D^{0,\,1}_\eta:C^\infty_{p,\,q}(X,\,L)\longrightarrow C^\infty_{p,\,q+1}(X,\,L)$ (one of which is used in (\ref{eqn:Laplacians_D_eta_0-1_star})) are defined. For example, for any $u\in C^\infty_{p,\,q}(X,\,L)$, we get \begin{eqnarray*}\overline{D^{0,\,1}_\eta}\,u=\overline{D^{0,\,1}_{\eta,\,\overline{L}}\,\overline{u}} = \overline{D_{\overline{L}}''\,\overline{u} - \frac{q}{\eta}\,\bar\partial\eta\wedge\overline{u}} = D'u - \frac{q}{\overline\eta}\,\partial\overline\eta\wedge u,\end{eqnarray*} where the first equality has already been mentioned above.

\begin{Conc}\label{Conc:D_0-1_eta_conjugate} The operator $\overline{D^{0,\,1}_\eta}:C^\infty_{p,\,q}(X,\,L)\longrightarrow C^\infty_{p+1,\,q}(X,\,L)$ is given by \begin{eqnarray*}\overline{D^{0,\,1}_\eta} = D' - \frac{q}{\overline\eta}\,\partial\overline\eta\wedge\cdot\,.\end{eqnarray*}

  Hence, its formal adjoint $\overline{D^{0,\,1}_\eta}^\star:C^\infty_{p+1,\,q}(X,\,L)\longrightarrow C^\infty_{p,\,q}(X,\,L)$ is given by \begin{eqnarray*}\overline{D^{0,\,1}_\eta}^\star = D'^\star - \frac{q}{\eta}\,(\partial\overline\eta\wedge\cdot\,.)^\star.\end{eqnarray*}

\end{Conc}

As for cohomology and harmonic forms, we only give the following statement for $D''$ and the operators it induces, although its analogue for $D'$ in place of $D''$, for $\overline{D^{0,\,1}_\eta}$ in place of $D^{0,\,1}_\eta$ and for $\Delta'_\eta$ in place of $\Delta''_\eta$ holds too and can be given the analogous proof. 

\begin{Obs}\label{Obs:D''_cohom_Hodge} Suppose that $D''^2 = 0$. Then:

\vspace{1ex} 

(i)\, $(D^{0,\,1}_\eta)^2 = 0$;

\vspace{1ex} 

(ii)\, for every bidegree $(p,\,q)$, the linear map \begin{eqnarray}\label{eqn:theta_eta_cohom}\theta_\eta:H^{p,\,q}_{D''}(X,\,L)\longrightarrow H^{p,\,q}_{D^{0,\,1}_\eta}(X,\,L), \hspace{5ex} \theta_\eta\bigg(\{u\}_{D''}\bigg):=\{\theta_\eta u\}_{D^{0,\,1}_\eta},\end{eqnarray} is {\bf well defined} and an {\bf isomorphism}, where the $D^{0,\,1}_\eta$-cohomology space and the $D''$-cohomology space are defined by \begin{eqnarray*}H^{p,\,q}_{D^{0,\,1}_\eta}(X,\,L) & := & \frac{\ker\bigg(D^{0,\,1}_\eta:C^\infty_{p,\,q}(X,\,L)\longrightarrow C^\infty_{p,\,q+1}(X,\,L)\bigg)}{\mbox{Im}\,\bigg(D^{0,\,1}_\eta:C^\infty_{p,\,q-1}(X,\,L)\longrightarrow C^\infty_{p,\,q}(X,\,L)\bigg)}\end{eqnarray*} and the analogous formula with $D''$ in place of $D^{0,\,1}_\eta$.

\vspace{1ex} 

(iii)\, if $X$ is {\bf compact}, for every bidegree $(p,\,q)$, we have the {\bf Hodge isomorphism} \begin{eqnarray}\label{eqn:Hodge_isom_Delta''_eta}H^{p,\,q}_{D^{0,\,1}_\eta}(X,\,L)\simeq{\cal H}_{\Delta''_\eta}^{p,\,q}(X,\,L):=\ker\bigg(\Delta''_\eta:C^\infty_{p,\,q}(X,\,L)\longrightarrow C^\infty_{p,\,q}(X,\,L)\bigg)\end{eqnarray} of {\bf finite-dimensional} vector spaces, mapping every $D^{0,\,1}_\eta$-cohomology class to its unique $\Delta''_\eta$-harmonic representative.

\end{Obs}

\noindent {\it Proof.} (i) is obtained either by a straightforward computation or as an immediate consequence of the formula $D^{0,\,1}_\eta = \theta_\eta D''\theta_\eta^{-1}$ (see (\ref{eqn:D_eta_0-1_def})).  

\vspace{1ex}

(ii) is obtained as in [Pop24, Proposition 2.3.] by noticing that \begin{eqnarray*}\theta_\eta(\ker D'') = \ker D_\eta^{0,\,1}  \hspace{5ex}\mbox{and}\hspace{5ex} \theta_\eta(\mbox{Im}\, D'') = \mbox{Im}\, D_\eta^{0,\,1}\end{eqnarray*} and by using the fact that $\theta_\eta$ is an isomorphism at the level of $L$-valued forms.

\vspace{1ex}

(iii) follows from standard harmonic theory via the ellipticity of $\Delta''_\eta$ and the compactness of $X$.  \hfill $\Box$

\section{Twisted commutation relations}\label{section:twisted-commutation_D-eta_0-1} The context is the same as in $\S$\ref{section:main-def}. The starting point is the following

\begin{Lem}\label{Lem:D-01_eta-star_tau-bar-star} The following formula holds for all $p,\,q$: \begin{eqnarray}\label{eqn::D-01_eta-star_tau-bar-star}(D^{0,\,1}_\eta)^\star + \overline\tau^\star = -i\,\bigg[\Lambda,\,D' + \frac{p}{\overline\eta}\,\partial\overline\eta\wedge\cdot\,\bigg]\end{eqnarray} on $C^\infty_{p,\,q+1}(X,\,L)$.

\end{Lem}    

\noindent {\it Proof.} Since $D''^\star = -i\,[\Lambda,\,D'] - \overline\tau^\star$ (cf. (ii) of (\ref{eqn:standard-comm-rel_bundle})) and $(\bar\partial\eta\wedge\cdot\,)^\star = i\,[\Lambda,\,\partial\overline\eta\wedge\cdot\,]$ (cf. (i) of Lemma \ref{Lem:com_2}), formula (\ref{eqn:D_eta_0-1_star}) translates to (\ref{eqn::D-01_eta-star_tau-bar-star}).  \hfill $\Box$

\vspace{2ex}

On the other hand, using Conclusion \ref{Conc:D_0-1_eta_conjugate}, for every $v\in C^\infty_{p,\,q+1}(X,\,L)$ we get: \begin{eqnarray*}\bigg[\Lambda,\,\overline{D^{0,\,1}_\eta}\bigg](v) & = & \Lambda\bigg(D'v - \frac{q+1}{\overline\eta}\,\partial\overline\eta\wedge v\bigg) - D'\Lambda(v) + \frac{q}{\overline\eta}\,\partial\overline\eta\wedge\Lambda(v).\end{eqnarray*} This means that \begin{eqnarray}\label{eqn:Lambda_D-0-1_eta-bar}-i\,\bigg[\Lambda,\,\overline{D^{0,\,1}_\eta}\bigg] = -i\,[\Lambda,\,D'] + \frac{q}{\overline\eta}\,\bigg[\Lambda,\,i\partial\overline\eta\wedge\cdot\,\bigg] + \frac{i}{\overline\eta}\,\Lambda(\partial\overline\eta\wedge\cdot\,) \hspace{3ex}\mbox{on}\hspace{1ex} C^\infty_{p,\,q+1}(X,\,L).\end{eqnarray}

Moreover, if we let, as in $\S4$ of [Pop24]: \begin{eqnarray*}\tau_\eta^{0,\,1}:=[\Lambda,\,D_\eta^{0,\,1}\omega\wedge\cdot\,] = \overline\tau - \bigg[\Lambda,\,\frac{1}{\eta}\,\bar\partial\eta\wedge\omega\wedge\cdot\,\bigg],\end{eqnarray*} where $D_\eta^{0,\,1}\omega: = \bar\partial\omega - \frac{1}{\eta}\,\bar\partial\eta\wedge\omega$ (so, we imply here that $D_\eta^{0,\,1}$ is defined by replacing $D''$ with $\bar\partial$ in (\ref{eqn:D_eta_0-1_def}) when it acts on scalar-valued forms), we get: \begin{eqnarray}\label{eqn:tau_eta_0-1_star} (\tau_\eta^{0,\,1})^\star = \overline\tau^\star - \bigg[\Lambda,\,\frac{1}{\eta}\,\bar\partial\eta\wedge\omega\wedge\cdot\,\bigg]^\star.\end{eqnarray} 

Putting together (\ref{eqn::D-01_eta-star_tau-bar-star}), (\ref{eqn:Lambda_D-0-1_eta-bar}) and (\ref{eqn:tau_eta_0-1_star}), we get the following identity on $C^\infty_{p,\,q+1}(X,\,L)$: \begin{eqnarray}\label{eqn:Lambda_D-0-1_eta-star_tau_eta_0-1_star}(D^{0,\,1}_\eta)^\star + (\tau_\eta^{0,\,1})^\star = -i\,\bigg[\Lambda,\,\overline{D^{0,\,1}_\eta}\bigg] - \frac{p+q}{\overline\eta}\,(\bar\partial\eta\wedge\cdot\,)^\star  - \bigg[\Lambda,\,\frac{1}{\eta}\,\bar\partial\eta\wedge\omega\wedge\cdot\,\bigg]^\star - \frac{i}{\overline\eta}\,\Lambda(\partial\overline\eta\wedge\cdot\,)\end{eqnarray} after using (i) of Lemma \ref{Lem:com_2} to infer that $\bigg[\Lambda,\,\bigg((p+q)/\overline\eta\bigg)\,i\partial\overline\eta\wedge\cdot\,\bigg] = \bigg((p+q)/\overline\eta\bigg)\,(\bar\partial\eta\wedge\cdot\,)^\star$.

Using this preliminary identity, we can now derive the first commutation relation in this context.

\begin{Prop}\label{Prop:eta-twisted_commutation-relations_1-0_0-1} Let $(X,\,\omega)$ be a complex Hermitian manifold with $\mbox{dim}_\C X=n\geq 2$ and let $(L,\,h,\,D=D' + D'')$ be a $C^\infty$ complex line bundle on $X$ equipped with a $C^\infty$ Hermitian fibre metric and with a linear connection.

  For any non-vanishing $C^\infty$ function $\eta:X\longrightarrow\C$, the following {\bf $\eta$-twisted commutation relations} hold on the $C^\infty$ $L$-valued differential forms of any (bi-)degree on $X$: \begin{eqnarray}\label{eqn:eta-twisted_commutation-relations_1-0_0-1}\nonumber & (a)& (D_\eta^{0,\,1})^\star + (\tau_\eta^{0,\,1})^\star = -i\,[\Lambda,\, \overline{D^{0,\,1}_\eta}] + \frac{n}{\overline\eta}\,\bigg[i\,\partial\overline\eta\wedge\cdot\,,\Lambda\bigg]; \\
  & (b) & \overline{D^{0,\,1}_\eta}^{\,\star} + \overline{\tau^{0,\,1}_\eta}^{\,\star} = i\,[\Lambda,\, D^{0,\,1}_\eta] - \frac{n}{\eta}\,\bigg[i\,\bar\partial\eta\wedge\cdot\,,\Lambda\bigg].\end{eqnarray}

\end{Prop}

\noindent {\it Proof.} It suffices to prove (a) since (b) follows from it by conjugation.

We fix an arbitrary bidegree $(p,\,q)$ and we prove (a) on $C^\infty_{p,\,q}(X,\,L)$. Let $k:=p+q$. Taking adjoints in the identity of Lemma \ref{Lem:com_3}, we get: $$\bigg[\Lambda,\,\frac{1}{\eta}\,\bar\partial\eta\wedge\omega\wedge\cdot\,\bigg]^\star = (n-k-1)\,\frac{1}{\overline\eta}\,(\bar\partial\eta\wedge\cdot\,)^\star - \frac{i}{\overline\eta}\,\partial\overline\eta\wedge\Lambda$$ on $(k+1)$-forms. Using this, (\ref{eqn:Lambda_D-0-1_eta-star_tau_eta_0-1_star}) translates to the following identities on $C^\infty_{p,\,q+1}(X,\,L)$: \begin{eqnarray*}(D^{0,\,1}_\eta)^\star + (\tau_\eta^{0,\,1})^\star & = & -i\,\bigg[\Lambda,\,\overline{D^{0,\,1}_\eta}\bigg] - \frac{p+q}{\overline\eta}\,(\bar\partial\eta\wedge\cdot\,)^\star  - \frac{i}{\overline\eta}\,\Lambda(\partial\overline\eta\wedge\cdot\,)  - \frac{n-k-1}{\overline\eta}\,(\bar\partial\eta\wedge\cdot\,)^\star + \frac{i}{\overline\eta}\,\partial\overline\eta\wedge\Lambda \\
  & = & -i\,\bigg[\Lambda,\,\overline{D^{0,\,1}_\eta}\bigg] - \frac{n-1}{\overline\eta}\,(\bar\partial\eta\wedge\cdot\,)^\star - \frac{1}{\overline\eta}\,[\Lambda,\,i\,\partial\overline\eta\wedge\cdot\,].\end{eqnarray*}

Since $[\Lambda,\,\partial\overline\eta\wedge\cdot\,] = -i\,(\bar\partial\eta\wedge\cdot\,)^\star$ (cf. (i) of Lemma \ref{Lem:com_2}), the above identity amounts to (a).  \hfill $\Box$

\vspace{2ex}

An immediate consequence of these $\eta$-twisted commutation relations is the following relation between the twisted Laplacians $\Delta''_\eta$ and $\Delta'_\eta$.

\begin{Prop}\label{Prop:eta-BKN_rough} Let $(X,\,\omega)$ be a complex Hermitian manifold with $\mbox{dim}_\C X=n\geq 2$ and let $(L,\,h,\,D=D' + D'')$ be a $C^\infty$ complex line bundle on $X$ equipped with a $C^\infty$ Hermitian fibre metric and with a linear connection.

  For any non-vanishing $C^\infty$ function $\eta:X\longrightarrow\C$, the following {\bf rough $\eta$-Bochner-Kodaira-Nakano ($\eta$-BKN) identity} holds on $C^\infty$ $L$-valued differential forms of any (bi-)degree on $X$: \begin{eqnarray}\label{eqn:eta-BKN_rough}\nonumber\Delta''_\eta = \Delta'_\eta + i\,\bigg[[D^{0,\,1}_\eta,\,\overline{D^{0,\,1}_\eta}],\,\Lambda\bigg] & + & \bigg[\overline{D^{0,\,1}_\eta},\,\overline{\tau^{0,\,1}_\eta}^{\,\star}\bigg] - \bigg[D^{0,\,1}_\eta,\,(\tau^{0,\,1}_\eta)^\star\bigg] \\
    & + & n\,\bigg[D^{0,\,1}_\eta,\,\frac{1}{\overline\eta}\,[i\partial\overline\eta\wedge\cdot\,,\,\Lambda]\bigg] + n\,\bigg[\overline{D^{0,\,1}_\eta},\,\frac{1}{\eta}\,[i\bar\partial\eta\wedge\cdot\,,\,\Lambda]\bigg].  \end{eqnarray}

\end{Prop}

\noindent {\it Proof.} Using the expression for $(D^{0,\,1}_\eta)^\star$ given in (a) of (\ref{eqn:eta-twisted_commutation-relations_1-0_0-1}), we get the second equality below: \begin{eqnarray}\label{eqn:eta-BKN_proof_1_bis}\Delta''_\eta & = & [D_\eta^{0,\,1},\,(D^{0,\,1}_\eta)^\star] = -i\,\bigg[D^{0,\,1}_\eta,\,[\Lambda,\,\overline{D^{0,\,1}_\eta}]\bigg] - [D^{0,\,1}_\eta,\,(\tau^{0,\,1}_\eta)^\star] + n\,\bigg[D_\eta^{0,\,1},\,\frac{1}{\overline\eta}\,[i\partial\overline\eta\wedge\cdot\,,\,\Lambda]\bigg].\end{eqnarray}

Now, the Jacobi identity yields the former equality below (expressing the first term on the r.h.s. of (\ref{eqn:eta-BKN_proof_1_bis})), while the $\eta$-twisted commutation relation (b) of (\ref{eqn:eta-twisted_commutation-relations_1-0_0-1}) yields the latter equality: \begin{eqnarray}\label{eqn:eta-BKN_proof_2_bis}\nonumber -i\,\bigg[D^{0,\,1}_\eta,\,[\Lambda,\,\overline{D^{0,\,1}_\eta}]\bigg] & = & \bigg[\overline{D^{0,\,1}_\eta},\,i\,[\Lambda,\,D^{0,\,1}_\eta]\bigg] + i\,\bigg[[\overline{D^{0,\,1}_\eta},\,D^{0,\,1}_\eta],\,\Lambda\bigg] \\
  & = & i\,\bigg[[D^{0,\,1}_\eta,\,\overline{D^{0,\,1}_\eta}],\,\Lambda\bigg] + \bigg[\overline{D^{0,\,1}_\eta},\,\overline{D^{0,\,1}_\eta}^\star + \overline{\tau^{0,\,1}_\eta}^\star + \frac{n}{\eta}\,[i\bar\partial\eta\wedge\cdot,\,\Lambda]\bigg].\end{eqnarray}

Plugging into (\ref{eqn:eta-BKN_proof_1_bis}) the expression given for $-i\,\bigg[D^{0,\,1}_\eta,\,[\Lambda,\,\overline{D^{0,\,1}_\eta}]\bigg]$ in (\ref{eqn:eta-BKN_proof_2_bis}) and using the equality $[\overline{D^{0,\,1}_\eta},\,\overline{D^{0,\,1}_\eta}^{\,\star}] = \Delta'_\eta$, we get (\ref{eqn:eta-BKN_rough}). \hfill $\Box$

\vspace{2ex}

It remains to compute the operator $i\,\bigg[[D^{0,\,1}_\eta,\,\overline{D^{0,\,1}_\eta}],\,\Lambda\bigg]$ that plays in (\ref{eqn:eta-BKN_rough}) the role of the curvature operator of the classical Bochner-Kodaira-Nakano identity.

\begin{Prop}\label{Prop:curvature_computation} Let $(X,\,\omega)$ be a complex Hermitian manifold with $\mbox{dim}_\C X=n\geq 2$ and let $(L,\,h,\,D=D' + D'')$ be a $C^\infty$ complex line bundle on $X$ equipped with a $C^\infty$ Hermitian fibre metric and with a linear connection. Fix an arbitrary non-vanishing $C^\infty$ function $\eta:X\longrightarrow\C$. 

 (a)\, For any bidegree $(p,\,q)$ and any $u^{p,\,q}\in C^\infty_{p,\,q}(X,\,L)$, the following identity holds: \begin{eqnarray}\label{eqn:curvature_computation_first}\nonumber i\,[D^{0,\,1}_\eta,\,\overline{D^{0,\,1}_\eta}]\,(u^{p,\,q}) = - \frac{1}{\eta}\,i\bar\partial\eta\wedge D' u^{p,\,q} - \frac{1}{\overline\eta}\,i\partial\overline\eta\wedge D'' u^{p,\,q} + \Omega_{p,\,q}\wedge u^{p,\,q},\end{eqnarray} where $\Omega_{p,\,q}$ is the smooth $\C$-valued $(1,\,1)$-form on $X$ defined by \begin{eqnarray*}\Omega_{p,\,q}:= i\Theta(D)^{1,\,1} + \frac{p-q}{|\eta|^2}\,i\partial\overline\eta\wedge\bar\partial\eta + \frac{p}{\eta^2}\,i\partial\eta\wedge\bar\partial\eta - \frac{q}{\overline\eta^2}\,i\partial\overline\eta\wedge\bar\partial\overline\eta - \frac{p}{\eta}\,i\partial\bar\partial\eta + \frac{q}{\overline\eta}\,i\partial\bar\partial\overline\eta\end{eqnarray*} and $\Theta(D)^{1,\,1}$ is the $(1,\,1)$-part of the {\bf curvature form} of the connection $D$, namely the $C^\infty$ scalar-valued $(1,\,1)$-form on $X$ defined by \begin{eqnarray*}D'D'' + D''D' = \Theta(D)^{1,\,1}\wedge\cdot\,.\end{eqnarray*}

  \vspace{1ex}

  (b)\, For any bidegree $(p,\,q)$ and any $u^{p,\,q}\in C^\infty_{p,\,q}(X,\,L)$, the following identity holds: \begin{eqnarray}\label{eqn:curvature_computation_second}\nonumber i\,\bigg[[D^{0,\,1}_\eta,\,\overline{D^{0,\,1}_\eta}],\,\Lambda\bigg]\,(u^{p,\,q}) & = & \frac{1}{\eta}\,[\Lambda,\,i\bar\partial\eta\wedge\cdot\,]\,(D' u^{p,\,q}) + \frac{1}{\overline\eta}\,[\Lambda,\,i\partial\overline\eta\wedge\cdot\,]\,(D'' u^{p,\,q}) \\
    \nonumber & + & \frac{1}{\eta}\,i\bar\partial\eta\wedge[\Lambda,\,D']\,(u^{p,\,q}) + \frac{1}{\overline\eta}\,i\partial\overline\eta\wedge[\Lambda,\,D'']\,(u^{p,\,q}) \\
 & + & [\Omega_{p,\,q}\wedge\cdot\,,\,\Lambda]\,(u^{p,\,q}).\end{eqnarray}

\end{Prop}

\noindent {\it Proof.} (a)\, We compute separately the two terms in the sum $[D^{0,\,1}_\eta,\,\overline{D^{0,\,1}_\eta}]\,(u^{p,\,q}) = (D^{0,\,1}_\eta\overline{D^{0,\,1}_\eta})\,(u^{p,\,q}) + (\overline{D^{0,\,1}_\eta}D^{0,\,1}_\eta)\,(u^{p,\,q})$.

\vspace{1ex}

$\bullet$ We use Conclusion \ref{Conc:D_0-1_eta_conjugate} to get the first line and then (\ref{eqn:D_eta_0-1_def}) to get the next line below: \begin{eqnarray*}(D^{0,\,1}_\eta\overline{D^{0,\,1}_\eta})\,(u^{p,\,q}) & = & D^{0,\,1}_\eta\bigg(D' u^{p,\,q} - \frac{q}{\overline\eta}\,\partial\overline\eta\wedge u^{p,\,q}\bigg) \\
  & = & D''\bigg(D' u^{p,\,q} - \frac{q}{\overline\eta}\,\partial\overline\eta\wedge u^{p,\,q}\bigg) - \frac{p+1}{\eta}\,\bar\partial\eta\wedge\bigg(D' u^{p,\,q} - \frac{q}{\overline\eta}\,\partial\overline\eta\wedge u^{p,\,q}\bigg).\end{eqnarray*}

This leads to \begin{eqnarray}\label{eqn:curvature_computation_first_proof_1}\nonumber(D^{0,\,1}_\eta\overline{D^{0,\,1}_\eta})\,(u^{p,\,q}) = D''D' u^{p,\,q} & + & \frac{q}{\overline\eta^2}\,\bar\partial\overline\eta\wedge\partial\overline\eta\wedge u^{p,\,q} - \frac{q}{\overline\eta}\,\bar\partial\partial\overline\eta\wedge u^{p,\,q} +  \frac{q}{\overline\eta}\,\partial\overline\eta\wedge D'' u^{p,\,q}\\
  & - & \frac{p+1}{\eta}\,\bar\partial\eta\wedge D'u^{p,\,q} + \frac{(p+1)q}{|\eta|^2}\,\bar\partial\eta\wedge\partial\overline\eta\wedge u^{p,\,q}.\end{eqnarray}

$\bullet$ Meanwhile, we use (\ref{eqn:D_eta_0-1_def}) to get the first line and then Conclusion \ref{Conc:D_0-1_eta_conjugate} to get the next line below: \begin{eqnarray*}(\overline{D^{0,\,1}_\eta}D^{0,\,1}_\eta)\,(u^{p,\,q}) & = & \overline{D^{0,\,1}_\eta}\bigg(D'' u^{p,\,q} - \frac{p}{\eta}\,\bar\partial\eta\wedge u^{p,\,q}\bigg) \\
  & = & D'\bigg(D'' u^{p,\,q} - \frac{p}{\eta}\,\bar\partial\eta\wedge u^{p,\,q}\bigg) - \frac{q+1}{\overline\eta}\,\partial\overline\eta\wedge\bigg(D'' u^{p,\,q} - \frac{p}{\eta}\,\bar\partial\eta\wedge u^{p,\,q}\bigg).\end{eqnarray*}

This leads to \begin{eqnarray}\label{eqn:curvature_computation_first_proof_2}\nonumber(\overline{D^{0,\,1}_\eta}D^{0,\,1}_\eta)\,(u^{p,\,q}) = D'D'' u^{p,\,q} & + & \frac{p}{\eta^2}\,\partial\eta\wedge\bar\partial\eta\wedge u^{p,\,q} - \frac{p}{\eta}\,\partial\bar\partial\eta\wedge u^{p,\,q} +  \frac{p}{\eta}\,\bar\partial\eta\wedge D'u^{p,\,q} \\
  & - & \frac{q+1}{\overline\eta}\,\partial\overline\eta\wedge D''u^{p,\,q} + \frac{p(q+1)}{|\eta|^2}\,\partial\overline\eta\wedge\bar\partial\eta\wedge u^{p,\,q}.\end{eqnarray}

\vspace{2ex}

Adding up (\ref{eqn:curvature_computation_first_proof_1}) and (\ref{eqn:curvature_computation_first_proof_2}) yields the equality claimed under (a).

\vspace{1ex}

(b)\, Starting from $i\,\bigg[[D^{0,\,1}_\eta,\,\overline{D^{0,\,1}_\eta}],\,\Lambda\bigg]\,(u^{p,\,q}) = i\,[D^{0,\,1}_\eta,\,\overline{D^{0,\,1}_\eta}]\,(\Lambda u^{p,\,q}) -\Lambda\bigg(i\,[D^{0,\,1}_\eta,\,\overline{D^{0,\,1}_\eta}]\,(u^{p,\,q})\bigg)$ and using the equality under (a), we get \begin{eqnarray}\label{eqn:curvature_computation_second_proof_1}\nonumber i\,\bigg[[D^{0,\,1}_\eta,\,\overline{D^{0,\,1}_\eta}],\,\Lambda\bigg]\,(u^{p,\,q}) & = & \bigg(-\frac{1}{\eta}\,i\bar\partial\eta\wedge D'\Lambda u^{p,\,q} + \frac{1}{\eta}\,\Lambda(i\bar\partial\eta\wedge D' u^{p,\,q})\bigg) \\
  \nonumber  & + & \bigg(-\frac{1}{\overline\eta}\,i\partial\overline\eta\wedge D''\Lambda u^{p,\,q} + \frac{1}{\overline\eta}\,\Lambda(i\partial\overline\eta\wedge D'' u^{p,\,q})\bigg) \\
  & + & \Omega_{p,\,q}\wedge\Lambda u^{p,\,q} - \Lambda(\Omega_{p,\,q}\wedge u^{p,\,q}).\end{eqnarray}

It remains to notice that the first parenthesis on the right of (\ref{eqn:curvature_computation_second_proof_1}) equals \begin{eqnarray*}\frac{1}{\eta}\,[\Lambda,\,i\bar\partial\eta\wedge\cdot\,]\,(D' u^{p,\,q}) + \frac{1}{\eta}\,i\bar\partial\eta\wedge[\Lambda,\,D']\,(u^{p,\,q}),\end{eqnarray*} while the second parenthesis on the right on the right of (\ref{eqn:curvature_computation_second_proof_1}) equals \begin{eqnarray*}\frac{1}{\overline\eta}\,[\Lambda,\,i\partial\overline\eta\wedge\cdot\,]\,(D'' u^{p,\,q}) + \frac{1}{\overline\eta}\,i\partial\overline\eta\wedge[\Lambda,\,D'']\,(u^{p,\,q}),\end{eqnarray*} to get (\ref{eqn:curvature_computation_second}).    \hfill $\Box$

\begin{Obs}\label{Obs:eta_real-valued_curvature-form} In the context of Proposition \ref{Prop:curvature_computation}, if $\eta$ is {\bf real-valued}, then \begin{eqnarray*}\Omega_{p,\,q}:= i\Theta(D)^{1,\,1} + (p-q)\,\gamma_\eta,\end{eqnarray*} where $\gamma_\eta$ is the $C^\infty$ real $(1,\,1)$-form on $X$ defined by \begin{eqnarray*}\gamma_\eta:=\frac{2}{\eta^2}\,i\partial\eta\wedge\bar\partial\eta -  \frac{1}{\eta}\,i\partial\bar\partial\eta.\end{eqnarray*}


\end{Obs}  

\noindent {\it Proof.} This is an immediate consequence of the definition of $\Omega_{p,\,q}$ given in (a) of Proposition \ref{Prop:curvature_computation}.  \hfill $\Box$


\section{Vanishing of certain line bundle cohomology groups}\label{section:vanishing_cohomology_bundle_compact} As a consequence of our $\eta$-Bochner-Kodaira-Nakano identity of Propositions \ref{Prop:eta-BKN_rough} and \ref{Prop:curvature_computation}, we get the following vanishing result. For the sake of generality, we give it for an arbitrary connection $D$, although in the (most interesting) case where $D$ is compatible with a given smooth Hermitian fibre metric $h$ on $L$ (i.e. $D$ is a {\it Hermitian} connection), our hypothesis $D''^2=0$ implies that $D'^2=0$ (as one can easily deduce from the property of the locally defined $1$-form associated with a Hermitian connection in a local trivialisation of $L$ given e.g. in [Dem97, V-$\S7.$]), hence also that the curvature form $\Theta(D)$ reduces to its $(1,\,1)$-component $\Theta(D)^{1,\,1}$. Thus, when $D$ is Hermitian, $\Theta(D)^{1,\,1}$ is closed and therefore $C_D=0$, a fact that simplifies hypothesis (\ref{eqn:_cohomology_bundle_compact_H1}) imposed on $\eta$. 

\begin{The}\label{The:vanishing_cohomology_bundle_compact} Let $X$ be a {\bf compact} complex manifold with $\mbox{dim}_\C X=n\geq 2$ and let $L$ be a $C^\infty$ complex line bundle on $X$ equipped with a linear connection $D=D' + D''$ such that $D''^2=0$. Fix integers $p,q\in\{0,\dots , n\}$ such that either $(p>q \hspace{1ex}\mbox{and}\hspace{1ex} p+q\geq n+1)$ or $(p<q \hspace{1ex}\mbox{and}\hspace{1ex} p+q\leq n-1)$.

  Suppose there exists a non-vanishing $C^\infty$ function $\eta:X\longrightarrow\R$ satisfying the following conditions: 

\vspace{1ex}

(i)\, the $C^\infty$ $(1,\,1)$-form $$\omega = \omega_{p,\,q,\,\eta}:= \frac{1}{p-q}\,i\Theta(D)^{1,\,1} + \gamma_\eta$$ is {\bf positive definite} at every point of $X$, where $\Theta(D)^{1,\,1}$ is the part of type $(1,\,1)$ of the curvature form $\Theta(D)$ of $D$ and $\gamma_\eta:=\frac{2}{\eta^2}\,i\partial\eta\wedge\bar\partial\eta - \frac{1}{\eta}\,i\partial\bar\partial\eta$;

\vspace{1ex}

(ii)\, the pointwise $\omega$-norm $|\partial\eta| = |\partial\eta|_\omega$ of the $(1,\,0)$-form $\partial\eta$ is {\bf small} relative to $|\eta|$ in that \begin{eqnarray}\label{eqn:_cohomology_bundle_compact_H1}C_1(\eta):=\sup\limits_X\frac{|\partial\eta|}{|\eta|} < \frac{1 - 4(n + \,\sqrt{n})\,C_D}{8(n + \,\sqrt{n})\,C(\varphi) + 14n + 10\,\sqrt{n} + 4n\,\sqrt{n}},\end{eqnarray} where $C_D:=\sup\limits_X|\partial(i\Theta(D)^{1,\,1})|_\omega$, $\varphi:=-\log|\eta|$ and $C(\varphi):=\sup\limits_X|i\partial\bar\partial\varphi|_\omega$.

\vspace{1ex}

Then, the $D''$-cohomology group of $L$ of bidegree $(p,\,q)$ {\bf vanishes}, namely $H^{p,\,q}_{D''}(X,\,L) = \{0\}$.

\end{The}

\noindent {\it Proof.} Fix an arbitrary $C^\infty$ Hermitian fibre metric $h$ on $L$. (No link is assumed between $h$ and $D$.) The positive definiteness assumption (i) amounts to $\omega$ defining a Hermitian metric on $X$. It is this metric $\omega$ that will be used (in conjunction with $h$ for $L$-valued forms) to define the pointwise inner product $\langle\,\cdot\,,\,\cdot\,\rangle$, the pointwise norm $|\,\cdot\,|$, the $L^2$-inner product $\langle\langle\,\cdot\,,\,\cdot\,\rangle\rangle$, the $L^2$-norm $||\,\cdot\,||$ and all the adjoints appearing in the expressions involved in this proof. Note that this special choice of metric $\omega$ means that $\omega$ is now closely related to the real $(1,\,1)$-form $\Omega_{p,\,q}:= i\Theta(D)^{1,\,1} + (p-q)\,\gamma_\eta$ of Observation \ref{Obs:eta_real-valued_curvature-form}. (Recall that in Propositions \ref{Prop:eta-BKN_rough} and \ref{Prop:curvature_computation} they were independent.)

We will exploit our $\eta$-BKN identity (\ref{eqn:eta-BKN_rough}) for the $\C$-line bundle $(L,\,h,\,D)\longrightarrow(X,\,\omega)$ to get a lower estimate for $\langle\langle\Delta''_\eta u^{p,\,q},\,u^{p,\,q}\rangle\rangle$ for all forms $u=u^{p,\,q}\in C^\infty_{p,\,q}(X,\,L)$. Fix such a form. 

\vspace{1ex}

$\bullet$ {\it Estimating the curvature term of (\ref{eqn:eta-BKN_rough}).}

\vspace{1ex}

Since $\Omega_{p,\,q} = (p-q)\,\omega$ and since $[\omega\wedge\cdot\,,\,\Lambda] = (p+q-n)\,\mbox{Id}$ on $(p,\,q)$-forms, (\ref{eqn:curvature_computation_second}) translates to \begin{eqnarray}\label{eqn:curvature-term_proof-ineq-BKN}i\,\bigg[[D^{0,\,1}_\eta,\,\overline{D^{0,\,1}_\eta}],\,\Lambda\bigg]\,(u^{p,\,q}) = A(u^{p,\,q}) + (p-q)\,(p+q-n)\,u^{p,\,q},\end{eqnarray} where the first-order operator $A$ is given (in our case of a real-valued $\eta$) on $C^\infty_{p,\,q}(X,\,L)$ by \begin{eqnarray*}\label{eqn:eta-BKN_1-0_0-1_curvature_2}A & = & \frac{1}{\eta}\,[\Lambda,\,i\bar\partial\eta\wedge\cdot\,]\,D' + \frac{1}{\eta}\,[\Lambda,\,i\partial\eta\wedge\cdot\,]\,D'' + \frac{1}{\eta}\,i\bar\partial\eta\wedge[\Lambda,\,D'] + \frac{1}{\eta}\,i\partial\eta\wedge[\Lambda,\,D''] \\
  & = & \frac{1}{\eta}\,[\Lambda,\,i\bar\partial\eta\wedge\cdot\,]\,\bigg(\overline{D_\eta^{0,\,1}} + \frac{q}{\eta}\,\partial\eta\wedge\cdot\,\bigg) + \frac{1}{\eta}\,[\Lambda,\,i\partial\eta\wedge\cdot\,]\,\bigg(D_\eta^{0,\,1} + \frac{p}{\eta}\,\bar\partial\eta\wedge\cdot\,\bigg) \\
  & - & \frac{1}{\eta}\,\bar\partial\eta\wedge\bigg(({D^{0,\,1}_\eta})^{\,\star} + \frac{p}{\eta}\,(\bar\partial\eta\wedge\cdot\,)^\star + \bar\tau^\star\bigg) + \frac{1}{\eta}\,\partial\eta\wedge\bigg(\overline{D^{0,\,1}_\eta}^{\,\star} + \frac{q}{\eta}\,(\partial\eta\wedge\cdot\,)^\star + \tau^\star\bigg).\end{eqnarray*} To get the latter equality above, we first used the commutation relations (ii) and (i) in (\ref{eqn:standard-comm-rel_bundle}) to express $[\Lambda,\,D']$ and $[\Lambda,\,D'']$ in terms of $D''^{\star}$ and $D'^{\star}$ and then (\ref{eqn:D_eta_0-1_star}) and Conclusion \ref{Conc:D_0-1_eta_conjugate} to further express $D''^{\star}$ and $D'^{\star}$ in terms of $({D^{0,\,1}_\eta})^{\,\star}$ and $\overline{D^{0,\,1}_\eta}^{\,\star}$.

Since \begin{eqnarray*}|\Lambda| = |L| = \sup_{|u|=1}|\omega\wedge u|\leq\sqrt{n}\end{eqnarray*} and \begin{eqnarray*}|\tau| =  \sup_{|u|=1}|\tau u| \leq \sup_{|u|=1}|\Lambda(\partial\omega\wedge u)| + \sup_{|u|=1}|\partial\omega\wedge\Lambda u| \leq 2\sqrt{n}\,|\partial\omega|\end{eqnarray*} (see [Pop24, proof of Theorem 5.2.]) and since in our case \begin{eqnarray*}\partial\omega = \frac{1}{p-q}\,\partial\bigg(i\Theta(D)^{1,\,1}\bigg)  + \partial\gamma_\eta = \frac{1}{p-q}\,\partial\bigg(i\Theta(D)^{1,\,1}\bigg) -\frac{2}{\eta^2}\,\partial\eta\wedge i\partial\bar\partial\eta,\end{eqnarray*} we get \begin{eqnarray}\label{eqn:tau_bound}|\tau| \leq 2\sqrt{n}\,\bigg(\frac{C_D}{|p-q|} + 2\,C_1(\eta)\,C_2(\eta)\bigg):=C_3(\eta)\end{eqnarray} at every point of $X$, where we put $C_2(\eta): = \sup\limits_X\frac{|\partial\bar\partial\eta|}{|\eta|}$.

Consequently, the expression obtained above for $A$ acting on $C^\infty_{p,\,q}(X,\,L)$ leads to the following inequalities for every $u = u^{p,\,q}\in C^\infty_{p,\,q}(X,\,L)$: \begin{eqnarray*}|\langle\langle Au,\,u\rangle\rangle| \leq \sqrt{n}\,C_1(\eta)\,||u||\, & \bigg(& ||\overline{D_\eta^{0,\,1}}u|| + q\,C_1(\eta)\,||u|| + ||D_\eta^{0,\,1}u|| + p\,C_1(\eta)\,||u|| + ||({D^{0,\,1}_\eta})^{\,\star} u|| \\
  & & + p\,C_1(\eta)\,||u|| + \,C_3(\eta)\,||u|| + ||\overline{D^{0,\,1}_\eta}^{\,\star}u|| + q\,C_1(\eta)\,||u|| + C_3(\eta)\,||u||\bigg).\end{eqnarray*} 

Applying the elementary inequality $ab\leq (a^2 + b^2)/2$ with $b=||u||$ and $a$ each of the $L^2$-norms of $\overline{D^{0,\,1}_\eta}u$, $D_\eta^{0,\,1}u$, $({D^{0,\,1}_\eta})^{\,\star} u$ and $\overline{D^{0,\,1}_\eta}^{\,\star}u$, we further get: \begin{eqnarray*}|\langle\langle Au,\,u\rangle\rangle| \leq \frac{\sqrt{n}}{2}\,C_1(\eta)\,\bigg(||\overline{D_\eta^{0,\,1}}u||^2 + ||\overline{D^{0,\,1}_\eta}^{\,\star}u||^2 + ||D_\eta^{0,\,1}u||^2 + ||({D^{0,\,1}_\eta})^{\,\star} u||^2 + C_{1,\,2}(p,\,q,\,n)\,||u||^2\bigg),\end{eqnarray*} where we put \begin{eqnarray*}C_{1,\,2}(p,\,q,\,n):= 4\,\bigg(1 + (p+q)\,C_1(\eta) + C_3(\eta)\bigg).\end{eqnarray*}

We have thus got: \begin{eqnarray}\label{eqn:vanishing_cohomology_bundle_compact_proof_1}|\langle\langle Au,\,u\rangle\rangle| \leq \frac{\sqrt{n}}{2}\,C_1(\eta)\,\bigg(\langle\langle\Delta'_\eta u,\,u\rangle\rangle + \langle\langle\Delta''_\eta u,\,u\rangle\rangle + C_{1,\,2}(p,\,q,\,n)\,||u||^2\bigg)\end{eqnarray} for every $u\in C^\infty_{p,\,q}(X,\,L)$.

\vspace{1ex}

$\bullet$ {\it Estimating the last term on the r.h.s. of (\ref{eqn:eta-BKN_rough})}.

\vspace{1ex}

Recall that $[i\bar\partial\eta\wedge\cdot\,,\,\Lambda] = (\partial\overline\eta\wedge\cdot\,)^\star$ (cf. (ii) of Lemma \ref{Lem:com_2}). Meanwhile, $\eta$ is real-valued in our case, so we get: \begin{eqnarray}\label{eqn:vanishing_cohomology_bundle_compact_proof_2}\nonumber \bigg|\bigg\langle\bigg\langle n\,\bigg[\overline{D^{0,\,1}_\eta},\,\frac{1}{\eta}\,[i\bar\partial\eta\wedge\cdot\,,\,\Lambda]\bigg] u,\,u\bigg\rangle\bigg\rangle\bigg| & \leq & n\,\bigg|\bigg\langle\bigg\langle\overline{D^{0,\,1}_\eta}\,u,\,\frac{1}{\eta}\,\partial\eta\wedge u\bigg\rangle\bigg\rangle\bigg| + n\,\bigg|\bigg\langle\bigg\langle\frac{1}{\eta}\,(\partial\eta\wedge\cdot\,)^\star u,\, \overline{D^{0,\,1}_\eta}^\star\,u\bigg\rangle\bigg\rangle\bigg| \\
\nonumber  & \leq & n\,C_1(\eta)\,\bigg(\bigg|\bigg|\overline{D^{0,\,1}_\eta}\,u\bigg|\bigg|\,||u|| + \bigg|\bigg|\overline{D^{0,\,1}_\eta}^\star\,u\bigg|\bigg|\,||u||\bigg) \\
\nonumber & \leq & \frac{n}{2}\,C_1(\eta)\,\bigg(\bigg|\bigg|\overline{D^{0,\,1}_\eta}\,u\bigg|\bigg|^2 + \bigg|\bigg|\overline{D^{0,\,1}_\eta}^\star\,u\bigg|\bigg|^2 + 2\,||u||^2\bigg) \\
 & = & \frac{n}{2}\,C_1(\eta)\,\bigg(\langle\langle\Delta'_\eta u,\,u\rangle\rangle + 2\,||u||^2\bigg)\end{eqnarray} for every $u\in C^\infty_{p,\,q}(X,\,L)$.

\vspace{1ex}

$\bullet$ {\it Estimating the last but one term on the r.h.s. of (\ref{eqn:eta-BKN_rough})}.

\vspace{1ex}

Recall that $[i\partial\overline\eta\wedge\cdot\,,\,\Lambda] = -(\bar\partial\eta\wedge\cdot\,)^\star$ (cf. (i) of Lemma \ref{Lem:com_2}). Thus, since $\overline\eta = \eta$, we get: \begin{eqnarray}\label{eqn:vanishing_cohomology_bundle_compact_proof_3}\nonumber\bigg|\bigg\langle\bigg\langle n\,\bigg[D^{0,\,1}_\eta,\,\frac{1}{\overline\eta}\,[i\partial\overline\eta\wedge\cdot\,,\,\Lambda]\bigg]\,u,\,u\bigg\rangle\bigg\rangle\bigg| & \leq & n\,\bigg|\bigg\langle\bigg\langle D^{0,\,1}_\eta\,u,\,\frac{1}{\eta}\,\bar\partial\eta\wedge u\bigg\rangle\bigg\rangle\bigg| + n\,\bigg|\bigg\langle\bigg\langle\frac{1}{\eta}\,(\bar\partial\eta\wedge\cdot\,)^\star u,\, (D^{0,\,1}_\eta)^\star\,u\bigg\rangle\bigg\rangle\bigg| \\
\nonumber  & \leq & n\,C_1(\eta)\,\bigg(\bigg|\bigg|D^{0,\,1}_\eta\,u\bigg|\bigg|\,||u|| + \bigg|\bigg|(D^{0,\,1}_\eta)^\star\,u\bigg|\bigg|\,||u||\bigg) \\
\nonumber & \leq & \frac{n}{2}\,C_1(\eta)\,\bigg(\bigg|\bigg|D^{0,\,1}_\eta\,u\bigg|\bigg|^2 + \bigg|\bigg|(D^{0,\,1}_\eta)^\star u\bigg|\bigg|^2 + 2\,||u||^2\bigg)  \\
 & = & \frac{n}{2}\,C_1(\eta)\,\bigg(\langle\langle\Delta''_\eta u,\,u\rangle\rangle + 2\,||u||^2\bigg)\end{eqnarray} for every $u\in C^\infty_{p,\,q}(X,\,L)$.

\vspace{1ex}

$\bullet$ {\it Estimating the two terms containing $\tau_\eta^{0,\,1}$ on the r.h.s. of (\ref{eqn:eta-BKN_rough})}.

\vspace{1ex}

Since $\tau_\eta^{0,\,1} = \overline\tau - \bigg[\Lambda,\,\frac{1}{\eta}\,\bar\partial\eta\wedge\omega\wedge\cdot\,\bigg]$, we get the next inequalities between pointwise operator norms: \begin{eqnarray*}|\tau_\eta^{0,\,1}|\leq|\tau| + \bigg|\bigg[\Lambda,\,\frac{1}{\eta}\,\bar\partial\eta\wedge\omega\wedge\cdot\,\bigg]\bigg| \leq C_3(\eta) + \bigg|\bigg[\Lambda,\,\frac{1}{\eta}\,\bar\partial\eta\wedge\omega\wedge\cdot\,\bigg]\bigg|.\end{eqnarray*}

Meanwhile, for every $u\in C^\infty_{p,\,q}(X,\,L)$, we have: \begin{eqnarray*}\bigg|\bigg\langle\bigg[\Lambda,\,\frac{1}{\eta}\,\bar\partial\eta\wedge\omega\wedge\cdot\,\bigg]\,u,\,u\bigg\rangle\bigg| & \leq & \bigg|\bigg\langle\frac{1}{\eta}\,\bar\partial\eta\wedge\omega\wedge u,\,\omega\wedge u\bigg\rangle\bigg| + \bigg|\bigg\langle\Lambda u,\,\bigg(\frac{1}{\eta}\,\bar\partial\eta\wedge\omega\wedge\cdot\,\bigg)^\star u\bigg\rangle\bigg| \\
  & \leq & n\,C_1(\eta)\,|u|^2 + n\,C_1(\eta)\,|u|^2.\end{eqnarray*}   

Therefore, we get the following upper bound for the pointwise operator norm of $\tau_\eta^{0,\,1}$: \begin{eqnarray*}|\tau_\eta^{0,\,1}|\leq C_3(\eta) + 2n\,C_1(\eta):=C_4(\eta).\end{eqnarray*} 

Using this, we get: \begin{eqnarray}\label{eqn:vanishing_cohomology_bundle_compact_proof_3}\nonumber\bigg|\bigg\langle\bigg\langle \bigg[D^{0,\,1}_\eta,\,(\tau^{0,\,1}_\eta)^\star\bigg]  \,u,\,u\bigg\rangle\bigg\rangle\bigg| & \leq & \bigg|\bigg\langle\bigg\langle D^{0,\,1}_\eta u,\,\tau^{0,\,1}_\eta u\bigg\rangle\bigg\rangle\bigg| + \bigg|\bigg\langle\bigg\langle(\tau^{0,\,1}_\eta)^\star u,\,(D^{0,\,1}_\eta)^\star u\bigg\rangle\bigg\rangle\bigg| \\
\nonumber & \leq & C_4(\eta)\,\bigg(\bigg|\bigg|D^{0,\,1}_\eta\,u\bigg|\bigg|\,||u|| + \bigg|\bigg|(D^{0,\,1}_\eta)^\star\,u\bigg|\bigg|\,||u||\bigg) \\
\nonumber & \leq & \frac{1}{2}\,C_4(\eta)\,\bigg(\bigg|\bigg|D^{0,\,1}_\eta\,u\bigg|\bigg|^2 + \bigg|\bigg|(D^{0,\,1}_\eta)^\star u\bigg|\bigg|^2 + 2\,||u||^2\bigg)  \\
& = & \frac{1}{2}\,C_4(\eta)\,\bigg(\langle\langle\Delta''_\eta u,\,u\rangle\rangle + 2\,||u||^2\bigg)\end{eqnarray} for every $u\in C^\infty_{p,\,q}(X,\,L)$.

Similarly, we get: \begin{eqnarray}\label{eqn:vanishing_cohomology_bundle_compact_proof_4}\nonumber\bigg|\bigg\langle\bigg\langle \bigg[\overline{D^{0,\,1}_\eta},\,\overline{\tau^{0,\,1}_\eta}^\star\bigg]  \,u,\,u\bigg\rangle\bigg\rangle\bigg| & \leq & \bigg|\bigg\langle\bigg\langle \overline{D^{0,\,1}_\eta} u,\,\overline{\tau^{0,\,1}_\eta} u\bigg\rangle\bigg\rangle\bigg| + \bigg|\bigg\langle\bigg\langle\overline{\tau^{0,\,1}_\eta}^\star u,\,\overline{D^{0,\,1}_\eta}^\star u\bigg\rangle\bigg\rangle\bigg| \\
\nonumber & \leq & C_4(\eta)\,\bigg(\bigg|\bigg|\overline{D^{0,\,1}_\eta}\,u\bigg|\bigg|\,||u|| + \bigg|\bigg|\overline{D^{0,\,1}_\eta}^\star\,u\bigg|\bigg|\,||u||\bigg) \\
\nonumber & \leq & \frac{1}{2}\,C_4(\eta)\,\bigg(\bigg|\bigg|\overline{D^{0,\,1}_\eta}\,u\bigg|\bigg|^2 + \bigg|\bigg|\overline{D^{0,\,1}_\eta}^\star u\bigg|\bigg|^2 + 2\,||u||^2\bigg)  \\
& = & \frac{1}{2}\,C_4(\eta)\,\bigg(\langle\langle\Delta'_\eta u,\,u\rangle\rangle + 2\,||u||^2\bigg)\end{eqnarray} for every $u\in C^\infty_{p,\,q}(X,\,L)$.

\vspace{1ex}

$\bullet$ Putting together (\ref{eqn:eta-BKN_rough}) and (\ref{eqn:curvature-term_proof-ineq-BKN})--(\ref{eqn:vanishing_cohomology_bundle_compact_proof_4}), we get: \begin{eqnarray*}\langle\langle\Delta''_\eta u,\,u\rangle\rangle & \geq & \langle\langle\Delta'_\eta u,\,u\rangle\rangle + (p-q)(p+q-n)\,||u||^2 \\
  & - & \frac{\sqrt{n}}{2}\,C_1(\eta)\,\bigg(\langle\langle\Delta'_\eta u,\,u\rangle\rangle + \langle\langle\Delta''_\eta u,\,u\rangle\rangle + C_{1,\,2}(p,\,q,\,n)\,||u||^2\bigg) \\
  & - & \frac{n}{2}\,C_1(\eta)\,\bigg(\langle\langle\Delta'_\eta u,\,u\rangle\rangle  + 2\,||u||^2\bigg) - \frac{n}{2}\,C_1(\eta)\,\bigg(\langle\langle\Delta''_\eta u,\,u\rangle\rangle  + 2\,||u||^2\bigg) \\
  & - & \frac{C_4(\eta)}{2}\,\bigg(\langle\langle\Delta''_\eta u,\,u\rangle\rangle  + 2\,||u||^2\bigg) - \frac{C_4(\eta)}{2}\,\bigg(\langle\langle\Delta'_\eta u,\,u\rangle\rangle  + 2\,||u||^2\bigg),\end{eqnarray*} which amounts to the following {\it twisted $\eta$-BKN inequality}: \begin{eqnarray}\label{eqn:eta-BKN_1-0_0-1_simplified_ineq}\bigg(1 + \frac{n+\sqrt{n}}{2}\,C_1(\eta) & + & \frac{C_4(\eta)}{2}\bigg)\,\langle\langle\Delta''_\eta u,\,u\rangle\rangle \geq \bigg(1 - \frac{n+\sqrt{n}}{2}\,C_1(\eta) - \frac{C_4(\eta)}{2}\bigg)\,\langle\langle\Delta'_\eta u,\,u\rangle\rangle \\
  \nonumber  & + & \bigg[(p-q)(p+q-n) - C_1(\eta)\,\bigg(2n + \frac{\sqrt{n}}{2}\,C_{1,\,2}(p,\,q,\,n)\bigg) - 2\,C_4(\eta)\bigg]\,||u||^2 \end{eqnarray} for every $u\in C^\infty_{p,\,q}(X,\,L)$.

\vspace{1ex}

$\bullet$ {\it Exploiting the twisted $\eta$-BKN inequality.}

Since $\langle\langle\Delta''_\eta u,\,u\rangle\rangle, \langle\langle\Delta'_\eta u,\,u\rangle\rangle\geq 0$, (\ref{eqn:eta-BKN_1-0_0-1_simplified_ineq}) shows that the vanishing of $\langle\langle\Delta''_\eta u,\,u\rangle\rangle$ (which is equivalent to $u\in\ker\Delta''_\eta$) implies the vanishing of $\langle\langle\Delta'_\eta u,\,u\rangle\rangle$ (which is irrelevant to us here) and the vanishing of $u$ whenever the coefficients of $\langle\langle\Delta'_\eta u,\,u\rangle\rangle$ and $||u||^2$ are positive. For this to happen, we need \begin{eqnarray*}\label{eqn:vanishing_cond_1} & (a) &  \frac{n+\sqrt{n}}{2}\,C_1(\eta) + \frac{C_4(\eta)}{2} < 1   \hspace{5ex}\mbox{and} \\
  & (b) & C_1(\eta)\,\bigg(2n + \frac{\sqrt{n}}{2}\,C_{1,\,2}(p,\,q,\,n)\bigg) + 2\,C_4(\eta) < (p-q)(p+q-n).\end{eqnarray*}

Since $p,q$ are integers and we require that either ($p>q$ and $p+q\geq n+1$) or ($p<q$ and $p+q\leq n-1$), the r.h.s. of (b) is a positive integer. Hence, (b) is satisfied whenever the inequality \begin{eqnarray*}(b')\hspace{2ex} C_1(\eta)\,\bigg(2n + \frac{\sqrt{n}}{2}\,C_{1,\,2}(p,\,q,\,n)\bigg) + 2\,C_4(\eta) < 1\end{eqnarray*} is satisfied. Meanwhile, $p+q\leq 2n$, so \begin{eqnarray*}C_{1,\,2}(p,\,q,\,n) \leq 4\bigg(1 + 2n\,C_1(\eta) + 4\sqrt{n}\,C_1(\eta)\,C_2(\eta) + 2\sqrt{n}\,C_D\bigg),\end{eqnarray*} while $1/|p-q|\leq 1$, so \begin{eqnarray*}C_4(\eta)\leq 2n\,C_1(\eta) + 4\sqrt{n}\,C_1(\eta)\,C_2(\eta) + 2\sqrt{n}\,C_D.\end{eqnarray*}

  Thus, for $(b')$ to hold, it suffices to have \begin{eqnarray*}C_1(\eta)\,\bigg[2n + 2\sqrt{n}\,\bigg(1 + 2n\,C_1(\eta) & + & 4\sqrt{n}\,C_1(\eta)\,C_2(\eta) + 2\sqrt{n}\,C_D\bigg)\bigg] \\
    & + & 2\,\bigg(2n\,C_1(\eta) + 4\sqrt{n}\,C_1(\eta)\,C_2(\eta) + 2\sqrt{n}\,C_D\bigg)<1.\end{eqnarray*} Hence, it suffices to have \begin{eqnarray*}C_1(\eta)^2\,\bigg(8n\,C_2(\eta) + 4n\sqrt{n}\bigg) + \bigg(8\sqrt{n}\,C_2(\eta) + 6n + 2\sqrt{n} + 4nC_D\bigg)\,C_1(\eta) + 4\sqrt{n}\,C_D  <1.\end{eqnarray*} Since $C_1(\eta)^2<C_1(\eta)$ (indeed, (a) implies $C_1(\eta)<1$), for this to happen it suffices that \begin{eqnarray*}(b'')\hspace{2ex} C_1(\eta)\,\bigg[8(n + \sqrt{n})\,C_2(\eta) + 6n + 4n\sqrt{n} + 2\sqrt{n} + 4n\,C_D\bigg] + 4\sqrt{n}\,C_D <1.\end{eqnarray*}

  Now, since $\varphi = -\frac{1}{2}\,\log(\eta^2)$, we have \begin{eqnarray}\label{eqn:eta-varphi}\partial\varphi = -\frac{1}{\eta}\,\partial\eta  \hspace{3ex}\mbox{and}\hspace{3ex} \frac{1}{\eta}\,\partial\bar\partial\eta = \partial\varphi\wedge\bar\partial\varphi - \partial\bar\partial\varphi.\end{eqnarray} Consequently, $C_1(\eta) = \sup\limits_X(|\partial\eta|/|\eta|) = \sup\limits_X|\partial\varphi|$ and \begin{eqnarray*}C_2(\eta) = \sup\limits_X\frac{|\partial\bar\partial\eta|}{\eta} \leq \sup\limits_X|\partial\varphi|^2 + \sup\limits_X|\partial\bar\partial\varphi| \leq C_1(\eta)^2 + C(\varphi).\end{eqnarray*} This upper estimate for $C_2(\eta)$ shows that $(b'')$ holds whenever \begin{eqnarray*} C_1(\eta)\,\bigg[8(n + \sqrt{n})\,\bigg(C_1(\eta)^2 + C(\varphi)\bigg) + 6n + 4n\sqrt{n} + 2\sqrt{n} + 4n\,C_D\bigg] + 4\sqrt{n}\,C_D <1.\end{eqnarray*}

  After expanding the l.h.s. of the above inequality and using the fact that $C_1(\eta) < 1$ (hence also $C_1(\eta)^3<C_1(\eta)$), we see that $(b'')$ (hence also (b)) holds whenever \begin{eqnarray*}(b''')\hspace{2ex} C_1(\eta)\,\bigg[8(n + \sqrt{n})\, C(\varphi) + 14n + 10\sqrt{n} + 4n\sqrt{n}\bigg] + 4(n+\sqrt{n})\,C_D <1.\end{eqnarray*} This condition is satisfied thanks to hypothesis (\ref{eqn:_cohomology_bundle_compact_H1}).

  To finish the proof of the implication \begin{eqnarray}\label{eqn:implication_Delta''-harmonic_zero}u\in\ker\Delta''_\eta\implies u=0,\end{eqnarray} it remains to observe that $(b''')$ implies (a). Indeed, (a) is implied by \begin{eqnarray*}(a')\hspace{2ex}\frac{3n+\sqrt{n}}{2}\,C_1(\eta) + 2\sqrt{n}\,C_1(\eta)\,C_2(\eta) + \sqrt{n}\,C_D  < 1.\end{eqnarray*} Since $C_2(\eta)\leq C_1(\eta)^2 + C(\varphi)$ and $C_1(\eta)^2<1$, we get $2\sqrt{n}\,C_1(\eta)\,C_2(\eta) < 2\sqrt{n}\,C_1(\eta)\,(1 + C(\varphi))$. Thus, $(a')$ is implied by \begin{eqnarray*}C_1(\eta)\,\bigg(\frac{3n+5\sqrt{n}}{2} + 2\sqrt{n}\,C(\varphi)\bigg) + \sqrt{n}\,C_D  < 1,\end{eqnarray*} which is implied by $(b''')$.

 \vspace{1ex}

 $\bullet$ {\it End of the proof.}

 Implication (\ref{eqn:implication_Delta''-harmonic_zero}) that was proved above means (in the notation of Observation \ref{Obs:D''_cohom_Hodge}) that \begin{eqnarray*}{\cal H}_{\Delta''_\eta}^{p,\,q}(X,\,L) = \{0\}.\end{eqnarray*} Thanks to (iii) of  Observation \ref{Obs:D''_cohom_Hodge}, this amounts to $H^{p,\,q}_{D^{0,\,1}_\eta}(X,\,L) = \{0\}$, which in turn amounts, thanks to the isomorphism (\ref{eqn:theta_eta_cohom}), to $H^{p,\,q}_{D''}(X,\,L) = \{0\}$.

 The proof of Theorem \ref{The:vanishing_cohomology_bundle_compact} is complete.  \hfill $\Box$

 \vspace{2ex}

 The following statement is a variant of Theorem \ref{The:vanishing_cohomology_bundle_compact} deduced thereof by making a slightly stronger assumption.

 \begin{Cor}\label{Cor:vanishing_cohomology_bundle_compact} In the setting of Theorem \ref{The:vanishing_cohomology_bundle_compact}, suppose there exists a non-vanishing $C^\infty$ function $\eta:X\longrightarrow\R$ such that the function $\varphi = -\log|\eta|$ satisfies the following conditions: 

\vspace{1ex}

(i)\, the $C^\infty$ $(1,\,1)$-form $$\rho = \rho_{p,\,q,\,\varphi}:= \frac{1}{p-q}\,i\Theta(D)^{1,\,1} + i\partial\bar\partial\varphi$$ is {\bf positive definite} at every point of $X$, where $\Theta(D)^{1,\,1}$ is the $(1,\,1)$-component of the curvature form of $D$;

\vspace{1ex}

(ii)\, the pointwise $\rho$-norm $|\partial\varphi| = |\partial\varphi|_\rho$ of the $(1,\,0)$-form $\partial\varphi$ is uniformly {\bf small} in that \begin{eqnarray}\label{eqn:_cor-cohomology_bundle_compact_H1}C_1'(\varphi):=\sup\limits_X|\partial\varphi|_\rho < \frac{1 - 4(n + \,\sqrt{n})\,C'_D}{8(n + \,\sqrt{n})\,C'(\varphi) + 14n + 10\,\sqrt{n} + 4n\,\sqrt{n}},\end{eqnarray} where $C'_D:=\sup\limits_X|\partial(i\Theta(D)^{1,\,1})|_\rho$ and $C'(\varphi):=\sup\limits_X|i\partial\bar\partial\varphi|_\rho$.

\vspace{1ex}

Then, the $D''$-cohomology group of $L$ of bidegree $(p,\,q)$ {\bf vanishes}, namely $H^{p,\,q}_{D''}(X,\,L) = \{0\}$.

\end{Cor}   

 \noindent {\it Proof.} In the notation of Theorem \ref{The:vanishing_cohomology_bundle_compact}, (\ref{eqn:eta-varphi}) gives the equality below: \begin{eqnarray*}\gamma_\eta = i\partial\varphi\wedge\bar\partial\varphi + i\partial\bar\partial\varphi \geq i\partial\bar\partial\varphi,\end{eqnarray*} the last inequality being a consequence of the standard inequality $i\alpha\wedge\overline\alpha\geq 0$ holding for any $(1,0)$-form $\alpha$. This implies (using again the notation of Theorem \ref{The:vanishing_cohomology_bundle_compact}) the inequality: \begin{eqnarray*}\omega:= \frac{1}{p-q}\,i\Theta(D)^{1,\,1} + \gamma_\eta \geq \frac{1}{p-q}\,i\Theta(D)^{1,\,1} + i\partial\bar\partial\varphi = \rho\end{eqnarray*} at every point of $X$.

   This inequality $\omega\geq\rho$ implies, on the one hand (thanks to our current hypothesis (i)): \begin{eqnarray*}\omega>0 \hspace{5ex} \mbox{on}\hspace{1ex} X\end{eqnarray*} and, on the other hand, the inequalities \begin{eqnarray*}C'(\varphi)\geq C(\varphi)  \hspace{5ex}\mbox{and}\hspace{5ex} C'_D\geq C_D,\end{eqnarray*} where $C(\varphi)$ and $C_D$ were defined in Theorem \ref{The:vanishing_cohomology_bundle_compact}. In particular, since our current hypothesis (ii) coincides with the second inequality below, we get: \begin{eqnarray*}C_1(\eta)\leq C_1'(\varphi) < \frac{1 - 4(n + \,\sqrt{n})\,C'_D}{8(n + \,\sqrt{n})\,C'(\varphi) + 14n + 10\,\sqrt{n} + 4n\,\sqrt{n}} \leq \frac{1 - 4(n + \,\sqrt{n})\,C_D}{8(n + \,\sqrt{n})\,C(\varphi) + 14n + 10\,\sqrt{n} + 4n\,\sqrt{n}},\end{eqnarray*} where $C_1(\eta)$ was defined in Theorem \ref{The:vanishing_cohomology_bundle_compact}.

     We conclude that the hypotheses of Theorem \ref{The:vanishing_cohomology_bundle_compact} are satisfied. Therefore, its conclusion, which coincides with our current conclusion, holds.  \hfill $\Box$

\vspace{2ex}

We now observe that Corollary \ref{Cor:vanishing_cohomology_bundle_compact} remains valid if we suppress its hypothesis (ii) but assume instead the connection $D$ to be {\it h-compatible}. This means that $D$ is a {\it Hermitian} connection, namely compatible with the Hermitian metric $h$ in the sense that \begin{eqnarray}\label{eqn:Hermitian-connection_def}d\{s,\,t\} = \{Ds,\,t\} + (-1)^p\,\{s,\,Dt\}\end{eqnarray} for every $p,q\in\{0,\dots, 2n\}$ and every $L$-valued $p$-form $s$ (i.e. every $s\in C^\infty_p(X,\,L)$) and every $L$-valued $q$-form $t$ (i.e. every $t\in C^\infty_q(X,\,L)$), where $\{\,\cdot\,,\,\cdot\,\}$ is the operation that combines the wedge product between the scalar parts of the $L$-valued differential forms with the pointwise inner product defined by $h$ between the vector parts of these forms. (See [Dem97, $\S V.7$].)   

\begin{Cor}\label{Cor:vanishing_cohomology_bundle_compact_simpler} In the setting of Corollary \ref{Cor:vanishing_cohomology_bundle_compact}, we assume furthermore that the connection $D$ is {\bf h-compatible}. Suppose there exists a $C^\infty$ function $\varphi:X\longrightarrow\R$ such that the $C^\infty$ $(1,\,1)$-form $$\rho = \rho_{p,\,q,\,\varphi}:= \frac{1}{p-q}\,i\Theta(D)^{1,\,1} + i\partial\bar\partial\varphi$$ is {\bf positive definite} at every point of $X$, where $\Theta(D)^{1,\,1}$ is the $(1,\,1)$-component of the curvature form $\Theta(D)$ of $D$.

 \vspace{1ex}

Then, the $D''$-cohomology group of $L$ of bidegree $(p,\,q)$ {\bf vanishes}, namely $H^{p,\,q}_{D''}(X,\,L) = \{0\}$.

\end{Cor}
  
\noindent {\it Proof.} We will observe that the hypotheses of Corollary \ref{Cor:vanishing_cohomology_bundle_compact} are satisfied after replacing $\varphi$ with $\varepsilon\varphi$ and $D$ with $\sqrt\varepsilon\,D$ for a sufficiently small positive constant $\varepsilon$.

Indeed, $\Theta(D)^{1,\,1}$ becomes $\varepsilon\,\Theta(D)^{1,\,1}$ (because $D'D'' + D''D' = \Theta(D)^{1,\,1}\wedge\cdot$ and $D'$, resp. $D''$, becomes $\sqrt\varepsilon\,D'$, resp. $\sqrt\varepsilon\,D''$), hence $\rho$ becomes $\varepsilon\rho$, which is still positive definite if $\rho$ is.

Meanwhile, $|\partial\varphi|_\rho$ becomes $$|\partial(\varepsilon\varphi)|_{\varepsilon\rho} = \frac{\varepsilon}{\sqrt\varepsilon}\,|\partial\varphi|_\rho = \sqrt\varepsilon\,|\partial\varphi|_\rho.$$ Therefore, $C_1'(\varphi)$ becomes $\sqrt\varepsilon\,C_1'(\varphi)$, which tends to $0$ when $\varepsilon\downarrow 0$.

On the other hand, $D$ is now supposed $h$-compatible. Consequently, the hypothesis $D''^2=0$ implies $D'^2=0$, hence $\Theta(D) = \Theta(D)^{1,\,1}$. Thus, $\Theta(D)^{1,\,1}$ is $d$-closed (because we always have $d\Theta(D) = 0$). This implies $\partial\Theta(D)^{1,\,1} = 0$ (since $\Theta(D)^{1,\,1}$ is of pure type), so $C'_D = 0$ (see the notation of Corollary \ref{Cor:vanishing_cohomology_bundle_compact}).

Finally, $C'(\varphi)$ remains unchanged since $\sup\limits_X|i\partial\bar\partial(\varepsilon\varphi)|_{\varepsilon\rho} = \frac{\varepsilon}{\varepsilon}\,\sup\limits_X|i\partial\bar\partial\varphi|_\rho = C'(\varphi)$.

We conclude that, when $\varepsilon\downarrow 0$, the r.h.s. of (\ref{eqn:_cor-cohomology_bundle_compact_H1}) remains constant, equal to \begin{eqnarray*}\frac{1}{8(n + \,\sqrt{n})\,C'(\varphi) + 14n + 10\,\sqrt{n} + 4n\,\sqrt{n}},\end{eqnarray*} while the l.h.s. of (\ref{eqn:_cor-cohomology_bundle_compact_H1}) equals $\sqrt\varepsilon\,C_1'(\varphi)\downarrow 0$ and can, therefore, be made arbitrarily small. Thus, hypothesis (\ref{eqn:_cor-cohomology_bundle_compact_H1}) is satisfied if $\varepsilon$ is chosen small enough.

Corollary \ref{Cor:vanishing_cohomology_bundle_compact} yields $H^{p,\,q}_{\sqrt\varepsilon\,D''}(X,\,L) = \{0\}$ for all small enough $\varepsilon$.

Since $H^{p,\,q}_{D''}(X,\,L) = H^{p,\,q}_{\sqrt\varepsilon\,D''}(X,\,L)$ for all $\varepsilon>0$, we are done.  \hfill $\Box$

\vspace{2ex}

\noindent {\it Proof of Theorem \ref{The:vanishing_cohomology_bundle_compact_introd}.} The hypotheses show that Theorem \ref{The:vanishing_cohomology_bundle_compact} can be applied to get Theorem \ref{The:vanishing_cohomology_bundle_compact_introd} by taking $\eta:=\eta_\varepsilon$ with $\varepsilon>0$ small enough.

Indeed, under the assumptions of Theorem \ref{The:vanishing_cohomology_bundle_compact_introd}, the l.h.s. of (\ref{eqn:_cohomology_bundle_compact_H1}) can be made arbitrarily close to $0$ by choosing $\eta:=\eta_\varepsilon$ and $\varepsilon$ small enough, while the r.h.s. of (\ref{eqn:_cohomology_bundle_compact_H1}) remains uniformly bounded below away from $0$: \begin{eqnarray*}\frac{1}{8(n + \,\sqrt{n})\,C(\varphi_\varepsilon) + 14n + 10\,\sqrt{n} + 4n\,\sqrt{n}} \geq \frac{1}{8(n + \,\sqrt{n})\,C + 14n + 10\,\sqrt{n} + 4n\,\sqrt{n}}, \hspace{5ex} \varepsilon>0.\end{eqnarray*} (As in Corollary \ref{Cor:vanishing_cohomology_bundle_compact}, the $h$-compatibility assumption on $D$ implies $C_D=0$.)

This means that hypothesis (ii) of Theorem \ref{The:vanishing_cohomology_bundle_compact} is satisfied when choosing $\eta:=\eta_\varepsilon$ and $\varepsilon$ small enough. The same goes for hypothesis (i) of Theorem \ref{The:vanishing_cohomology_bundle_compact} thanks to the analogous assumption (i) of Theorem \ref{The:vanishing_cohomology_bundle_compact_introd}. \hfill $\Box$

\vspace{2ex}

Finally, we observe that the proof of Theorem \ref{The:vanishing_cohomology_bundle_compact} gives information on the lower part of the spectrum of the twisted Laplacian $\Delta''_\eta$ in certain bidegrees. This ties in with the discussion initiated in [Pop24] for scalar-valued forms (see Problem 1.1. and Proposition 1.3. therein).

\begin{Cor}\label{Cor:spectrum-lower-bound} The hypotheses are those of Theorem \ref{The:vanishing_cohomology_bundle_compact}. The following statements hold.

\vspace{1ex}

(i)\, The spectrum of $\Delta''_\eta$ in bidegree $(p,\,q)$ has the lower bound given by the inclusion: \begin{eqnarray}\label{eqn:spectrum-lower-bound_p-q}\mbox{Spec}\,\bigg(\Delta''_\eta : C^\infty_{p,\,q}(X,\,L)\longrightarrow C^\infty_{p,\,q}(X,\,L)\bigg)\subset\bigg[B_{p,\,q,\,n},\,+\infty\bigg),\end{eqnarray} where  \begin{eqnarray*}B_{p,\,q;\,n}:=\frac{(p-q)(p+q-n) - C_1(\eta)\,\bigg(2n + \frac{\sqrt{n}}{2}\,C_{1,\,2}(p,\,q,\,n)\bigg) - 2\,C_4(\eta)}{1 + \frac{n+\sqrt{n}}{2}\,C_1(\eta) + \frac{C_4(\eta)}{2}}\end{eqnarray*} and $C_1(\eta)$ is defined in (\ref{eqn:_cohomology_bundle_compact_H1}), $C_{1,\,2}(p,\,q,\,n) = 4\,\bigg(1 + (p+q)\,C_1(\eta) + C_3(\eta)\bigg)$, $C_4(\eta) = C_3(\eta) + 2n\,C_1(\eta)$, while $C_3(\eta)$ is defined in (\ref{eqn:tau_bound}).

\vspace{1ex}

(ii)\, The positive part of the spectrum of $\Delta''_\eta$ in bidegree $(n,\,0)$ has the lower bound given by the inclusion: \begin{eqnarray}\label{eqn:spectrum-lower-bound_n-0}\mbox{Spec}\,\bigg(\Delta''_\eta : C^\infty_{n,\,0}(X,\,L)\longrightarrow C^\infty_{n,\,0}(X,\,L)\bigg)\subset\{0\}\cup\bigg[B_{n,\,1;\,n},\,+\infty\bigg),\end{eqnarray} where $B_{n,\,1;\,n}$ is the constant defined under (i) for $(p,\,q) = (n,\,1)$.

\end{Cor}  

\noindent {\it Proof.} Part (i) is an immediate consequence of the twisted $\eta$-BKN inequality (\ref{eqn:eta-BKN_1-0_0-1_simplified_ineq}) in which all the coefficients are positive thanks to properties (a) and (b) spelt out in the proof of Theorem \ref{The:vanishing_cohomology_bundle_compact} being satisfied (as shown there, as a consequence of the hypotheses).

To prove (ii), we first use Observation \ref{Obs:D''_cohom_Hodge} to infer, from the hypothesis $D''^2 = 0$, that $(D^{0,\,1}_\eta)^2 = 0$. This implies that $D^{0,\,1}_\eta$ commutes with $\Delta''_\eta$. Hence, whenever $u\in C^\infty_{p,\,q}(X,\,L)\setminus\{0\}$ such that $\Delta''_\eta u = \lambda u$ for some real number $\lambda$, we have $$\Delta''_\eta(D^{0,\,1}_\eta u) = \lambda\, D^{0,\,1}_\eta u.$$

This implies that, if $\lambda\neq 0$ and $(p,\,q) = (n,\,0)$ (or merely $q=0$ and $p$ is arbitrary), $D^{0,\,1}_\eta u$ is a $\lambda$-eigenvector for $\Delta''_\eta$ in bidegree $(n,\,1)$ whenever $u$ is a $\lambda$-eigenvector for $\Delta''_\eta$ in bidegree $(n,\,0)$. Indeed, to prove this it remains to show that $D^{0,\,1}_\eta u \neq 0$ whenever $u\neq 0$. Suppose that $D^{0,\,1}_\eta u = 0$. Since $u$ is of bidegree $(n,\,0)$, we would also have (trivially, for bidegree reasons) $(D^{0,\,1}_\eta)^\star u = 0$. Therefore, we would have $\Delta''_\eta u = 0$. Thanks to the choice of $u$, this would mean that $\lambda u = 0$, which is impossible when both $\lambda$ and $u$ are non-zero.

Thus, every non-zero eigenvalue $\lambda$ of $\Delta''_\eta : C^\infty_{n,\,0}(X,\,L)\longrightarrow C^\infty_{n,\,0}(X,\,L)$ is again an eigenvalue for $\Delta''_\eta : C^\infty_{n,\,1}(X,\,L)\longrightarrow C^\infty_{n,\,1}(X,\,L)$. Along with (\ref{eqn:spectrum-lower-bound_p-q}) for $(p,\,q) = (n,\,1)$, this proves (\ref{eqn:spectrum-lower-bound_n-0}).  \hfill $\Box$

\section{Appendix: review of standard commutation relations}\label{Appendix}

 We briefly recall here some standard formulae that were used throughout the paper.

 \begin{Lem}\label{Lem:com_1} Let $(X,\,\omega)$ be a complex Hermitian manifold.

(a)\, The following standard Hermitian commutation relations ([Dem84], see also [Dem97, VI, $\S6.2$]) hold: \begin{eqnarray}\label{eqn:standard-comm-rel}\nonumber &  & (i)\,\,(\partial + \tau)^{\star} = i\,[\Lambda,\,\bar\partial];  \hspace{3ex} (ii)\,\,(\bar\partial + \bar\tau)^{\star} = - i\,[\Lambda,\,\partial]; \\
&  & (iii)\,\, \partial + \tau = -i\,[\bar\partial^{\star},\,L]; \hspace{3ex} (iv)\,\,
\bar\partial + \bar\tau = i\,[\partial^{\star},\,L],\end{eqnarray}

   \noindent where the upper symbol $\star$ stands for the formal adjoint w.r.t. the $L^2$ inner product induced by $\omega$, $L=L_{\omega}:=\omega\wedge\cdot$ is the Lefschetz operator of multiplication by $\omega$, $\Lambda=\Lambda_{\omega}:=L^{\star}$ and $\tau:=[\Lambda,\,\partial\omega\wedge\cdot]$ is the torsion operator (of order zero and type $(1,\,0)$) associated with the metric $\omega$.

\vspace{1ex}

(b)\, Let $(E,\,h)$ be a Hermitian $C^\infty$ $\C$-vector bundle of arbitrary rank over $X$. Let $D=D' + D''$ be an arbitrary connection on $E$, where $D'$ and $D''$ are its components of types $(1,\,0)$, respectively $(0,\,1)$. The operators $L$ and $\Lambda$ are extended to $E$-valued forms in $\Lambda^{p,\,q}T^\star X\otimes E$ by taking their tensor products with $\mbox{Id}_E$, while $\tau$ continues to be defined by $\tau:=[\Lambda,\,\partial\omega\wedge\cdot]$.

The Hermitian commutation relations of (a) still hold in the $E$-valued context (see [Dem97, VII, $\S.1$]), namely we have: \begin{eqnarray}\label{eqn:standard-comm-rel_bundle}\nonumber &  & (i)\,\,(D' + \tau)^{\star} = i\,[\Lambda,\,D''];  \hspace{3ex} (ii)\,\,(D'' + \bar\tau)^{\star} = - i\,[\Lambda,\,D']; \\
&  & (iii)\,\, D' + \tau = -i\,[D''^{\star},\,L]; \hspace{3ex} (iv)\,\,
D'' + \bar\tau = i\,[D'^{\star},\,L].\end{eqnarray} 

\end{Lem}

\vspace{2ex}

The following formulae can be viewed as commutation relations for zeroth-order operators. They are the generalisations for complex-valued functions of the formulae given in [Pop03, $\S.1.0.2$] and [Pop23, Appendix] for real-valued functions. 

\begin{Lem}\label{Lem:com_2} Let $(X,\,\omega)$ be a complex Hermitian manifold and $\eta$ a {\bf complex-valued} ${\cal C}^\infty$ function on $X$. The following identities hold pointwise for arbitrary differential forms of any degree on $X$. \\

\hspace{6ex} (i)\, $[\partial\bar\eta\wedge\cdot, \, \Lambda]=i\, (\bar\partial\eta\wedge\cdot)^{\star}$; \hspace{4ex}  (ii)\, $[\bar\partial\eta\wedge\cdot, \, \Lambda] = -i\, (\partial\bar\eta\wedge\cdot)^{\star}$; \\

\hspace{6ex} (iii)\, $[L, \, (\partial\bar\eta\wedge\cdot)^{\star}] = -i\, \bar\partial\eta\wedge\cdot$;  \hspace{4ex} (iv)\, $[L, \,(\bar\partial\eta\wedge\cdot)^{\star}] = i\, \partial\bar\eta\wedge\cdot$. 


\end{Lem}

\noindent {\it Proof.} It suffices to prove (i) as the other formulae follow from it by taking conjugates and adjoints.

All the operators involved in these identities are defined pointwise, so it suffices to prove them at an arbitrary given point. Moreover, we will show that the identities hold even when the adjoints are taken w.r.t. the pointwise inner product $\langle\,\cdot\,,\,\cdot\,\rangle$ induced by $\omega$.

Let $x\in X$ be an arbitrary point and let $z_1, \dots , z_n$ be local holomorphic coordinates about $x$ such that $\omega(x) = \sum_k i\, dz_k \wedge d\bar{z}_k$. Then, $(d\bar{z}_j\wedge\cdot)^\star = (\partial/\partial\bar{z}_j)\lrcorner\cdot$ at $x$ w.r.t. the pointwise inner product induced by $\omega$ at $x$. For any $(p,\,q)$-form $u$ and any $(p,\,q+1)$-form $v$, we get at $x$: \begin{eqnarray*}\langle u,\, (\bar\partial\eta\wedge\cdot\,)^\star v\rangle & = & \langle \bar\partial\eta\wedge u, \,v\rangle = \sum\limits_{j=1}^n\bigg\langle\frac{\partial\eta}{\partial\bar{z}_j}\,d\bar{z}_j\wedge u, \,v\bigg\rangle = \bigg\langle u, \, \sum\limits_{j=1}^n \frac{\partial\bar\eta}{\partial z_j}\,\bigg(\frac{\partial}{\partial\bar{z}_j}\lrcorner v\bigg)\bigg\rangle.\end{eqnarray*} This proves the first of the equivalent, mutually conjugated, formulae below: \begin{eqnarray}\label{eqn:com_2_proof_1}(\bar\partial\eta\wedge\cdot\,)^\star = \sum\limits_{j=1}^n \frac{\partial\bar\eta}{\partial z_j}\,\bigg(\frac{\partial}{\partial\bar{z}_j}\lrcorner\cdot\,\bigg) \iff (\partial\bar\eta\wedge\cdot\,)^\star = \sum\limits_{j=1}^n \frac{\partial\eta}{\partial\bar{z}_j}\,\bigg(\frac{\partial}{\partial z_j}\lrcorner\cdot\,\bigg).\end{eqnarray} 

On the other hand, for every form $v$, we have at $x$: \begin{eqnarray}\label{eqn:com_2_proof_2}\nonumber\bigg(i\,[\Lambda,\,\partial\bar\eta\wedge\cdot\,]\bigg)^\star v & = & -i\,[(\partial\bar\eta\wedge\cdot\,)^\star,\,\omega\wedge\cdot\,]\,(v) = -i\,(\partial\bar\eta\wedge\cdot\,)^\star(\omega\wedge v) + i\,\omega\wedge(\partial\bar\eta\wedge\cdot\,)^\star v  \\
  & = & -i\,\sum\limits_{j=1}^n\frac{\partial\eta}{\partial\bar{z}_j}\,\bigg(\frac{\partial}{\partial z_j}\lrcorner(\omega\wedge v)\bigg) + i\,\omega\wedge\sum\limits_{j=1}^n\frac{\partial\eta}{\partial\bar{z}_j}\,\bigg(\frac{\partial}{\partial z_j}\lrcorner v\bigg),\end{eqnarray} where the last equality followed from the previous one and from (\ref{eqn:com_2_proof_1}).

Now, still at $x$, we have: \begin{eqnarray*}\frac{\partial}{\partial z_j}\lrcorner(\omega\wedge v) = \bigg(\frac{\partial}{\partial z_j}\lrcorner\omega\bigg)\wedge v + \omega\wedge\bigg(\frac{\partial}{\partial z_j}\lrcorner v\bigg) = i\,d\bar{z}_j\wedge v + \omega\wedge\bigg(\frac{\partial}{\partial z_j}\lrcorner v\bigg),   \hspace{3ex} j=1,\dots , n.\end{eqnarray*} Using this, (\ref{eqn:com_2_proof_2}) becomes: \begin{eqnarray*}\bigg(i\,[\Lambda,\,\partial\bar\eta\wedge\cdot\,]\bigg)^\star v = \bigg(\sum\limits_{j=1}^n\frac{\partial\eta}{\partial\bar{z}_j}\,d\bar{z}_j\bigg)\wedge v = \bar\partial\eta\wedge v.\end{eqnarray*}

We have thus proved that $i\,[\Lambda,\,\partial\bar\eta\wedge\cdot\,] = (\bar\partial\eta\wedge\cdot\,)^\star$, which amounts to the stated identity (i). \hfill $\Box$





\vspace{2ex}

The following formula is a complement to the above ones.

\begin{Lem}\label{Lem:com_3} Let $(X,\,\omega)$ be a complex Hermitian manifold with $\mbox{dim}_\C X = n$ and $\eta$ a non-vanishing {\bf complex-valued} ${\cal C}^\infty$ function on $X$. Fix any $k\in\{0,\dots , n\}$.

  The following identity holds pointwise for arbitrary differential forms of degree $k$ on $X$: $$\bigg[\Lambda,\,\frac{1}{\eta}\,\bar\partial\eta\wedge\omega\wedge\cdot\,\bigg]= (n-k-1)\,\frac{1}{\eta}\,\bar\partial\eta\wedge\cdot\, + \omega\wedge\bigg(\frac{i}{\eta}\,(\partial\overline\eta\wedge\cdot\,)^\star\bigg).$$

\end{Lem}

\noindent {\it Proof.} Let $u$ be a $k$-form on $X$. We have: \begin{eqnarray*}\bigg[\Lambda,\,\frac{1}{\eta}\,\bar\partial\eta\wedge\omega\wedge\cdot\,\bigg]\,(u) & = & \Lambda\bigg(\omega\wedge\frac{1}{\eta}\,\bar\partial\eta\wedge u\bigg) - \omega\wedge\bigg(\frac{1}{\eta}\,\bar\partial\eta\wedge\Lambda(u)\bigg) \\
  & = & [\Lambda,\,\omega\wedge\cdot\,]\,\bigg(\frac{1}{\eta}\,\bar\partial\eta\wedge u\bigg) + \omega\wedge\bigg[\Lambda,\,\frac{1}{\eta}\,\bar\partial\eta\wedge\cdot\,\bigg]\,(u) \\
 & = & (n-k-1)\,\bigg(\frac{1}{\eta}\,\bar\partial\eta\wedge u\bigg) + \omega\wedge\bigg(\frac{i}{\eta}\,(\partial\overline\eta\wedge\cdot\,)^\star\bigg)\,(u),\end{eqnarray*} where the last equality follows from the standard identity $[\Lambda,\,\omega\wedge\cdot\,] = (n-k)\,\mbox{Id}$ on $k$-forms (applied here to the $(k+1)$-form $(1/\eta)\,\bar\partial\eta\wedge u$) and from (ii) of Lemma \ref{Lem:com_2}.  \hfill $\Box$

\vspace{6ex}

\noindent {\bf References} \\

\noindent [AN54]\, Y. Akizuki, S. Nakano --- {\it Note on Kodaira-Spencer's Proof of Lefschetz Theorems} --- Proc. Jap. Acad., {\bf 30} (1954), 266-272. 

\vspace{1ex}

\noindent [Dem84]\, J.-P. Demailly --- {\it Sur l'identit\'e de Bochner-Kodaira-Nakano en g\'eom\'etrie hermitienne} --- S\'eminaire d'analyse P. Lelong, P. Dolbeault, H. Skoda (editors) 1983/1984, Lecture Notes in Math., no. {\bf 1198}, Springer Verlag (1986), 88-97.

\vspace{1ex}

\noindent [Dem85]\, J.-P. Demailly --- {\it Champs magn\'etiques et in\'egalit\'es de Morse pour la $d''$-cohomologie} --- Ann. Inst. Fourier {\bf 35}, no. 4 (1985), 189-229.

\vspace{1ex}

\noindent [Dem 97]\, J.-P. Demailly --- {\it Complex Analytic and Algebraic Geometry}---http://www-fourier.ujf-grenoble.fr/~demailly/books.html

\vspace{1ex}

\noindent [DP04]\, J.-P. Demailly, M. Paun --- {\it Numerical Charaterization of the K\"ahler Cone of a Compact K\"ahler Manifold} --- Ann. of Math. (2) {\bf 159(3)} (2004) 1247-1274.



\vspace{1ex}

\noindent [GR70]\, H. Grauert, O. Riemenschneider --- {\it Verschwindungss\"atze f\"ur analytische Kohomologiegruppen auf komplexen R\"aumen} --- Invent. Math. {\bf 11} (1970), 263-292.

\vspace{1ex}

\noindent [Kod53]\, K. Kodaira --- {\it On a Differential Geometric Method in the Theory of Analytic Stacks} --- Proc. Nat. Acad. USA {\bf 39} (1953), 1268-1273.

\vspace{1ex}

\noindent [LP75]\, J. Le Potier --- {\it Annulation de la cohomologie \`a valeurs dans un fibr\'e vectoriel holomorphe de rang quelconque} --- Math. Ann., {\bf 218} (1975), 35-53.

\vspace{1ex}

\noindent [Nak55]\, S. Nakano --- {\it On Complex Analytic Vector Bundles} --- J. Math. Soc. Japan, {\bf 7} (1955), 1-12.

\vspace{1ex}

\noindent [Pop03]\, D. Popovici --- {\it Quelques applications des m\'ethodes effectives en g\'eom\'etrie analytique} --- PhD thesis, University Joseph Fourier (Grenoble 1), http://tel.ccsd.cnrs.fr/documents/

\vspace{1ex}

\noindent [Pop19]\, D. Popovici --- {\it Adiabatic Limit and the Fr\"olicher Spectral Sequence} --- Pacific Journal of Mathematics, {\bf 300}, No. 1, 2019, dx.doi.org/10.2140/pjm.2019.300.121

\vspace{1ex}

\noindent [Pop23]\, D. Popovici --- {\it A Non-Integrable Ohsawa-Takegoshi-Type $L^2$ Extension Theorem} --- arXiv e-print CV 2309.11291v1.

\vspace{1ex}

\noindent [Pop24]\, D. Popovici --- {\it Twisted Adiabatic Limit for Complex Structures} --- arXiv: 2404.06908v1 [math.DG]

\vspace{1ex}

\noindent [Siu84]\, Y.-T. Siu --- {\it A Vanishing Theorem for Semipositive Line Bundles over Non-K\"ahler Manifolds} --- J. Diff. Geom. {\bf 19} (1984) 431-452.

\vspace{1ex}

\noindent [Siu85]\, Y.-T. Siu --- {\it Some Recent Results in Complex Manifold Theory Related to Vanishing Theorems for the Semipositive Case} --- Springer Lecture Notes in Math {\bf 1111}, 169-192 (1985).

\vspace{6ex}

\noindent Universit\'e Paul Sabatier, Institut de Math\'ematiques de Toulouse

\noindent 118, route de Narbonne, 31062, Toulouse Cedex 9, France

\noindent Email: popovici@math.univ-toulouse.fr

\end{document}